\RequirePackage{amsmath} 
\documentclass[12pt]{iopart}

\usepackage{iopams}  
\usepackage{bbm}
\usepackage{amsthm}
\usepackage{amssymb}
\usepackage{tikz}
\usepackage{pgfplots}
\usepackage{cleveref}

\renewcommand{\phi}{\varphi}

\Crefname{equation}{}{}

\newtheorem{theorem}{Theorem}
\newtheorem{lemma}[theorem]{Lemma}
\newtheorem{corollary}[theorem]{Corollary}

\newtheorem{remark}[theorem]{Remark}
\newtheorem{proposition}[theorem]{Proposition}
\newtheorem{conj}{Conjecture}

\begin{document}

\title{On trajectories of complex-valued interior transmission eigenvalues}

\author{Lukas Pieronek$^1$ and Andreas Kleefeld$^{2,3}$}

\address{$^1$ Karlsruhe Institute of Technology, Institute for Applied and Numerical Mathematics, Englerstr. 2, 76131 Karlsruhe, Germany}
\address{$^2$ Forschungszentrum J\"ulich GmbH, 
	J\"ulich Supercomputing Centre, Wilhelm-Johnen-Str., 52425 J\"ulich, Germany}
\address{$^3$ University of Applied Sciences Aachen, Faculty of Medical Engineering and Technomathematics, Heinrich-Mu\ss{}mann-Str. 1, 52428 J\"ulich, Germany}

\eads{\mailto{l.pieronek@fz-juelich.de} and \mailto{a.kleefeld@fz-juelich.de}}
\vspace{10pt}
\begin{indented}
\item[]October 2022
\end{indented}

\begin{abstract}
	This paper investigates the interior transmission problem for homogeneous media via eigenvalue trajectories parameterized by the magnitude of the refractive index. In case that the scatterer is the unit disk, we prove that there is a one-to-one correspondence between complex-valued interior transmission eigenvalue trajectories and Dirichlet eigenvalues of the Laplacian which turn out to be exactly the trajectorial limit points as the refractive index tends to infinity. For general simply-connected scatterers in two or three dimensions, a corresponding relation is still open, but further theoretical results and numerical studies indicate a similar connection.
	\\[0.5cm]
\end{abstract}



\section{Introduction}
While Dirichlet eigenvalues of the Laplacian (DELs) are among the most famous and long-understood eigenvalues in PDE history, the first appearance of interior transmission eigenvalues (ITEs) reaches back to 1986, cf. \cite{kirsch1986denseness}, when Kirsch studied denseness properties of the far field operator in the context of inverse scattering problems.
Accordingly, ITEs correspond to critical and scatterer-specific wave numbers for which the feasibility of many shape-reconstructing sampling methods cannot be ensured. Apart from its physical origins, the related eigenproblem \Cref{ITP} --- the interior transmission problem (ITP) --- has also attracted own interests from a functional analytical perspective due to its non-selfadjoint nature. ITEs therefore require, in comparison with DELs that are classified by a self-adjoint operator, quite a non-standard approach of mathematical investigation which is why they also exhibit surprising structural phenomena. One of those will be the focus of this article and addresses a link between non-real ITEs and DELs of the same scattering object.

The importance of such a connection arises from the fact that only little is known about complex-valued ITEs in general. One of the very first results in this direction for quite arbitrary scatterers $D$ was the evidence of discreteness of the ITE spectrum and that infinity is its only accumulation point, see \cite{rynne1991interior,colton1989far,nguyen2017discreteness}. Further, it is known that for sufficiently smooth scatterers and indexes of refraction $n$ with $n\ne 1$ on $\partial D$ all ITEs are located in a horizontal strip around the real axis, see \cite{vodev2017high}, and there is an ITE-free lemniscate region centered at zero, see \cite{cakoni2010interior}. For the special case of a disk in 2D or a ball in 3D as scatterer, existence of non-real ITEs was proven for spherically-stratified media in \cite{leung2012complex}, providing thereby a more detailed distribution analysis of ITEs, see \cite{colton2015distribution}, and also first results for the inverse spectral problem, see \cite{colton2013complex}. However, existence of non-real ITEs for scatterers other than the disk is still an open problem. 
Still, there is numerical evidence which indicates existence of non-real interior transmission eigenvalues for arbitrary scatterers, see for example \cite[Table 3]{gintides2013computational}, \cite[Section 6.10]{kleefeld2013numerical}, \cite[Tables 4 and 5]{geng2016c}, \cite[Section 5]{yang2016mixed}, \cite[Section 5]{yang2016non}, \cite[Section 4]{han2016adaptive}, \cite[Tables 1 and 2]{han2015h}, \cite[Section 4]{xi2017recursive}, \cite[Table 5]{blaasten2017vanishing}, \cite[Tables 1 and 2]{han2017new}, \cite[Tables 3 and 4]{yang2017c0ipg}, \cite[Section 5]{li2018adaptive}, \cite[Section 5]{wang2018}, and \cite[Section 5]{xijizhang} to mention just a few.

The current paper takes a completely novel approach to the ITP and considers trajectories of complex-valued ITEs for homogeneous media which are parametrized by the magnitude of the refractive index $n$.
We show in case that $D$ is the unit disk that any non-real ITE generates a continuous trajectory when varying $n$ which is also well-defined as $n$ tends to infinity. Our main result is that there is a one-to-one correspondence between complex-valued ITE trajectories and DELs of the unit disk in the sense that for any DEL there exists exactly one complex-conjugated pair of complex-valued ITE trajectories which converge to that DEL as $n\to\infty$. This correspondence even holds in terms of geometric eigenvalue multiplicity, that is, the eigenspace dimension of any DEL coincides with the number of overlapping ITE trajectories from linear independent ITP eigenfunctions. 
As a byproduct, we obtain an alternative proof for the existence of complex-valued ITEs of the unit disk for constant $n$.
For quite arbitrary scatterers $D$, we prove that the only possible accumulation points of complex-valued ITE trajectories as $n\to \infty$ are also DELs, restricting thus the possible localisation of non-real ITEs for large refractive indices. 
Our numerical results for some simply-connected scatterers additionally indicate that a one-to-one correspondence between complex-valued ITE trajectories and DELs including multiplicity as $n\to \infty$ might even hold in more generality. These observations together finally give rise to an intrinsic relation between non-real ITEs and DELs which we conclude with in a conjecture.

The rest of this paper is organized as follows: after a short introduction to the interior transmission problem, \Cref{sec2} presents our theoretical findings about complex ITE trajectories with a strong focus on the unit disk. \Cref{sec3} provides complementing numerical results based on which we then formulate a conjecture on the observed relation between complex-valued ITE trajectories and DELs for simply-connected scatterers.

\section{Theoretical results}\label{sec2}

The interior transmission problem in the acoustic regime reads 
\begin{align}
	\begin{split}\label{ITP}
		\Delta w + \kappa^2n w=0 \quad &\text{in }D\ ,\\
		\Delta v + \kappa^2 v=0 \quad &\text{in }D\ ,\\
		v=w \quad &\text{on }\partial D\ , \\
		\partial_{\nu}v=\partial_{\nu}w \quad &\text{on }\partial D\ ,
	\end{split}
\end{align}
and is a system of Helmholtz equations coupled through the boundary data which impose non-scattering transmission conditions for time-harmonic waves. Here, $D$ is a bounded domain, which we also refer to as scatterer, and $n\in L^\infty(D)$ denotes the index of refraction that we assume to be constant and positive throughout this work. We call a wave number $\kappa=\kappa_n\in \mathbb{C}\backslash\{0\}$ for $n\ne 1$ an ITE if there exist non-trivial $v_n,w_n\in L^2(D)$ solving \Cref{ITP} such that $(v_n-w_n)\in H^2_0(D)$. 
For the special case of the unit disk $D\subset \mathbb{R}^2$ which most of this section will be devoted to ITP eigenfunction pairs are spanned by Fourier Bessel functions for $r\in [0,1)$ and $\phi\in[0,2\pi)$ given by
\begin{align}\label{eigf}
	\begin{split}
		v_n(r,\phi)&=J_p(\kappa_nr)\cos(p\phi)\ , \quad \Big(v_n(r,\phi)=J_p(\kappa_nr)\sin(p\phi)\ , p> 0\Big)\ , \\
		w_n(r,\phi)&=\alpha_n J_p(\sqrt{n}\kappa_nr)\cos(p\phi)\ , \quad \Big(w_n(r,\phi)=\alpha_nJ_p(\sqrt{n}\kappa_nr)\mathrm\sin(p\phi)\ , p> 0\Big)\ ,
	\end{split}
\end{align}
where $p\in \mathbb{N}_0$ is the Bessel order, $\alpha_n\in \mathbb{C}\backslash \{0\}$ is a coefficient to match the ITP boundary conditions in \Cref{ITP} and $J_p$ are Bessel functions of the first kind solving the Bessel equation
\begin{align}\label{bessel}
	z^2J''_p(z)+zJ'_p(z)+(z^2-p^2)J_p(z)=0\ .
\end{align}
Hence, ITEs $\kappa_n$ of the unit disk are equivalently characterized as roots of
\begin{align}
	F_p(n,\kappa):=\kappa \left(J_p'(\kappa)J_p(\kappa \sqrt{n})- \sqrt{n}J_p(\kappa)J_p'(\kappa \sqrt{n})\right)\ .
	\label{determinant}
\end{align} 	
To put ourselves into the framework of eigenvalue trajectories for the unit disk as scatterer, recall from the implicit function theorem that in a local open neighborhood of any $n$ subject to $\partial_\kappa F_p(n,\kappa)\ne 0$ there is a unique continuously-differentiable mapping $n\mapsto \kappa_n$ (which we will also denote by $\kappa_n$ and its derivative with respect to $n$ by $\kappa'_n$) such that $F_p(n,\kappa)=0$ if and only if $\kappa=\kappa_n$. 
A simple calculation shows that
\begin{align}\label{delk}
	\partial_\kappa F_p(n,\kappa)=(n-1)\kappa J_p(\kappa)J_p(\sqrt{n}\kappa)
\end{align}
and using the fact that $J_p$ only has real-valued roots we can already infer that ITE trajectories are indeed locally well-defined around non-real ITEs of the unit disk. For the rest of this paper, a continuous  ITE trajectory $\kappa_n$ is refered to as being complex-valued in some interval $(n_{\min},n_{\max})$ if $\kappa_n\in\mathbb{C}\backslash\mathbb{R}$ for $n\in(n_{\min},n_{\max})$ almost everywhere. We start with a qualitative characteristic of non-real ITEs for quite general scatterers $D$. 

\begin{lemma}\label{zerolemma}
	Let $D$ be a bounded domain with $C^{1,1}$-boundary and  $\kappa_n\in\mathbb{C}\backslash \mathbb{R}$ be an ITE with eigenfunctions $v_n,w_n\in L^2(D)$. Then it holds that
	\begin{align}\label{vanish}
		\int_{D}\vert v_n\vert^2-n\vert w_n\vert^2\, \mathrm{d}x= 0\ .
	\end{align}
\end{lemma}
\begin{proof}
	The result is taken from \cite{Pi20} with a slightly different proof. We set $u_n:=(v_n-w_n)\in H^2_0(D)$ and compute 
	\begin{align*}
		&\ \kappa_{n}^2\int_{D} v_{n}\overline{v}_{n} - n\overline{w}_{n}w_{n}\, \mathrm{d}x \\
		=&\int_{D} -\Delta v_{n}\overline{v}_{n} +  \overline{w}_{n}\Delta w_{n}\, \mathrm{d}x \\
		=&\int_{D} -\Delta v_{n}\overline{v}_{n} +  \overline{w}_{n}\Delta w_{n}\, \mathrm{d}x 
		+\int_{D} \Delta w_{n}\overline{v}_{n} -  \overline{v}_{n}\Delta w_{n}\, \mathrm{d}x\\
		=&\int_{D}-\Delta u_{n} \overline{v}_{n} -  \overline{u}_{n}\Delta w_{n}\, \mathrm{d}x  \\
		=&\int_{D}-u_{n} \Delta\overline{v}_{n} -  \Delta\overline{u}_{n}  w_{n}\, \mathrm{d}x  \\
		=&\ \overline{\kappa}_{n}^2\int_{D} -(v_n-w_n)(-\overline{v}_{n})-(-\overline{v}_{n} + n \overline{w}_{n})w_{n}\, \mathrm{d}x \\
		=&\ \overline{\kappa}_{n}^2\int_{D} v_{n}\overline{v}_{n} - n\overline{w}_{n}w_{n}\, \mathrm{d}x\ .
	\end{align*}
	Note that the boundary terms from integration by parts, which is well-defined for distributional solutions $v_n,w_n$ of the Helmholtz equation since $D$ is sufficiently regular,
	vanish due to $u_{n}\in H^2_0(D)$.
	Since $\kappa_n\ne \overline{\kappa}_n$ by assumption, we conclude that 
	\begin{align*}
		\int_{D}\vert v_n\vert^2-n\vert w_n\vert^2\, \mathrm{d}x= \int_{D} v_{n}\overline{v}_{n} - n\overline{w}_{n}w_{n}\, \mathrm{d}x= 0\ .
	\end{align*}
\end{proof}

\noindent Complementarily, we specify \eqref{vanish} for real-valued ITEs of the unit disk.

\begin{lemma}\label{energy}
	Let $\kappa_n$ be an ITE of the unit disk $D$ for $n\ne 1$ whose eigenfunction pair $(v_n,w_n)$ is given by \Cref{eigf} for some $p\in \mathbb{N}_0$. Then it holds that
	\begin{align}\label{difference}
		\int_{D}\vert v_n\vert^2-n\vert w_n\vert^2\, \mathrm{d}x=
		\left\{
		\begin{array}{ll} 
			\frac{(1-n)}{2}J_p(\kappa_n)^2 \int_0^{2\pi}\cos(p\phi)^2\,\mathrm{d}\phi\ , & \text{for }\kappa_n\in \mathbb{R}\ ,\\
			0\ , & \text{else}\ .
		\end{array}
		\right.
	\end{align}
\end{lemma}
\begin{proof}
	Let $(v_n,w_n)$ be as in \Cref{eigf} with cosine as angular part. If $\kappa_n\not\in \mathbb{R}$, it follows by the previous lemma that 
	\begin{align*}
		\int_{D}|v_n|^2-n|w_n|^2\,\mathrm{d}x=0\ .
	\end{align*}
	In the other case, that is $\kappa_n\in \mathbb{R}$, we obtain by using polar coordinates
	\begin{align*}
		\int_{D}v_{n}^2-n w_{n}^2\,\mathrm{d}x
		=\int_0^{2\pi}\cos(p\phi)^2\,\mathrm{d}\phi \int_0^1\left(J_p(\kappa_{n}r)^2-n\alpha_{n}^2 J_p(\sqrt{n}\kappa_{n}r)^2\right)r\,\mathrm{d}r\ .
	\end{align*}
	Next, we employ formula (5.14.5) from \cite[p. 129]{Le70} to evaluate
	\begin{align}\label{int1}
		\int_0^1 J_p(\kappa_{n}r)^2r\,\mathrm{d}r=\frac{1}{2}\left[J'_p(\kappa_{n})^2+\left(1-\frac{p}{\kappa_{n}^2}\right)J_p(\kappa_{n})^2\right]
	\end{align}
	and similarly, including the ITP boundary conditions via $\alpha_n^2\ne 0$,
	\begin{align}\label{int2}
		\begin{split}
			n\int_0^1 \alpha_{n}^2J_p(\sqrt{n}\kappa_{n}r)^2r\,\mathrm{d}r
			&=\frac{1}{2}\left[n\alpha_{n}^2J'_p(\sqrt{n}\kappa_{n})^2+n\left(1-\frac{p}{n\kappa_{n}^2}\right)\alpha_{n}^2J_p(\sqrt{n}\kappa_{n})^2\right]\\
			&=\frac{1}{2}\left[J'_p(\kappa_{n})^2+\left(n-\frac{p}{\kappa_{n}^2}\right)J_p(\kappa_{n})^2\right]\ .
		\end{split}
	\end{align}
	Subtracting \Cref{int2} from \Cref{int1} yields
	\begin{displaymath}
		\int_{D}v_{n}^2-n w_{n}^2\,\mathrm{d}x=\frac{(1-n)}{2}J_p(\kappa_{n})^2 \int_0^{2\pi}\cos(p\phi)^2\,\mathrm{d}\phi\ .
	\end{displaymath}
	This completes the proof if cosine is the angular part of the eigenfunction pair. The case of sine with $p>0$ follows along the same lines, noting that
	\begin{displaymath}
		\int_0^{2\pi}\cos(p\phi)^2\,\mathrm{d}\phi=\int_0^{2\pi}\sin(p\phi)^2\,\mathrm{d}\phi=\pi\ .
	\end{displaymath}
\end{proof}
\noindent \Cref{energy} yields an immediate conclusion for intersection points of complex-valued ITE trajectories of the unit disk with the real axis.
\begin{corollary}\label{cor}
	Let $\kappa_n$ be a continuous ITE trajectory of the unit disk $D$ for $n\in (n^\ast-\epsilon,n^\ast+\epsilon)$ and $\epsilon>0$ sufficiently small such that $1\not\in (n^\ast-\epsilon,n^\ast+\epsilon)$ whose eigenfunction pairs $(v_n,w_n)$ are given by \Cref{eigf} for some fixed $p\in \mathbb{N}_0$. If $\kappa_{n^\ast}\in \mathbb{R}$, but $\kappa_n\in \mathbb{C}\backslash \mathbb{R}$ for all $n\in(n^\ast-\epsilon,n^\ast+\epsilon)\backslash \{n^\ast\}$, then $\kappa_{n^\ast}$ is a DEL of the unit disk with eigenfunction $v_{n^\ast}$ and $\sqrt{n^\ast}\kappa_{n^\ast}$ is a DEL of the unit disk with eigenfunction $w_{n^\ast}$.
\end{corollary}
\begin{proof}
	The assertion that $\kappa_{n^\ast}$ is a DEL with eigenfunction $v_{n^\ast}$ follows by \Cref{difference} and continuity of $n\mapsto \int_{D}|v_n|^2-n|w_n|^2\,\mathrm{d}x$ which enforces $J_p(\kappa_{n^\ast})=0$. Comparing with the ITP boundary conditions \Cref{ITP} then, also $\sqrt{n^\ast}\kappa_{n^\ast}$ is a DEL with  eigenfunction $w_{n^\ast}$ since $\alpha_n\ne 0$ for all $n\ne 1$.
\end{proof}

\noindent The following lemma elaborates on the local behavior of ITE trajectories of the unit disk near DELs and shows that complex-valued trajectories hitting DELs is a frequent event accompanied by additional real-valued trajectories each.

\begin{lemma}\label{big}
	Let $\kappa_n$ be a continuous and complex-valued ITE trajectory of the unit disk $D$ for $n\in (n^\ast-\epsilon,n^\ast+\epsilon)$ with $\epsilon>0$ sufficiently small such that $1\not\in (n^\ast-\epsilon,n^\ast+\epsilon)$ whose eigenfunction pairs $(v_n,w_n)$ are given by \Cref{eigf} for some fixed $p\in \mathbb{N}_0$. If $J_p(\kappa_{n^\ast})=0$, then there exists a further continuous ITE trajectory $\widetilde{\kappa}_n$ defined in a local neighborhood of $n^\ast$ which is real-valued, assigned to the same $p\in \mathbb{N}_0$ and also fulfills $\widetilde{\kappa}_{n^\ast}=\kappa_{n^\ast}$. 
	Conversely, if  $\kappa^\ast$ is such that $J_p(\kappa^\ast)=0$, then there exist infinitely many $n^\ast\nearrow \infty$, a complex-conjugated pair of complex-valued ITE trajectories $\kappa_n, \overline{\kappa}_n$ as well as a real-valued ITE trajectory $\widetilde{\kappa}_n$ all of which are assigned to the same $p\in \mathbb{N}_0$, defined in a local neighborhood of $n^\ast$, are continuous and fulfill $\kappa_{n^\ast}=\widetilde{\kappa}_{n^\ast}=\kappa^\ast$.
\end{lemma}
\begin{proof}
	For the first assertion, recall from \Cref{delk} that
	\begin{align}\label{deriv}
		\partial_\kappa F_p(n,\kappa)=(n-1)\kappa J_p(\kappa)J_p(\sqrt{n}\kappa)\ ,
	\end{align}
	which vanishes at $(n^\ast,\kappa_{n^\ast})$ since $J_p(\kappa_{n^\ast})=J_p(\sqrt{n}\kappa_{n^\ast})=0$ by assumption. Using $J_p'(\kappa_{n^\ast})\ne 0 \ne J_p'(\sqrt{n}\kappa_{n^\ast})$ we deduce that the multiplicity of the root $\kappa=\kappa_{n^\ast}$ within the holomorphic function $\kappa\mapsto \partial_\kappa F_p(n^\ast,\kappa)$ is 2 and therefore 3 for $\kappa\mapsto F_p(n^\ast,\kappa)$. By Hurwitz' theorem we know that the total number of zeros of $\kappa\mapsto F_p(n,\kappa)$, including multiplicity, remains 3 in a small neighborhood of $\kappa=\kappa_{n^\ast}$ for all $n$ sufficiently close to $n^\ast$. Since $\kappa_n$ is complex-valued for $n\ne n^\ast$ by assumption, it approaches $\kappa_{n^\ast}$ as $n\to n^\ast$ in complex-conjugated pairs, so their combined number of roots is even. Hence,  there must be some real-valued trajectory $\widetilde{\kappa}_n$ which also goes to $\kappa_{n^\ast}$ as $n\to n^\ast$. Continuity of $\widetilde{\kappa}_n$ also follows by Hurwitz theorem for $n=n^\ast$ and for $n\ne n^\ast$ by the implicit function theorem. 
	
	Conversely, let $\kappa^\ast$ be a DEL of the unit disk such that $J_p(\kappa^\ast)=0$. Since $J_p$ has infinitely many roots, we can find $\kappa^{\ast\ast}\ne \kappa^\ast$ such that $J_p(\kappa^{\ast\ast})=0$ and set
	\begin{align}\label{ratio}
		n^\ast:=(\kappa^{\ast\ast}/\kappa^\ast)^2\ .
	\end{align}
	By the same reasoning as above, the total number of roots of $\kappa \mapsto F_p(n^\ast,\kappa)$ including multiplicity remains 3 in a small neighborhood of $\kappa^\ast$ for $n$ sufficiently close to $n^\ast$. Thus there exists trajectories $\kappa_{n,1},\kappa_{n,2},\kappa_{n,3}$ for $n\in (n^\ast-\epsilon,n^\ast+\epsilon)$ for some $\epsilon>0$ and associated with Bessel index $p$ such that $\kappa_{n^\ast,1}=\kappa_{n^\ast,2}=\kappa_{n^\ast,3}=\kappa^\ast$. These are, in particular, continuous in $n^\ast$ and continuously differentiable for $n\ne n^\ast$ according to the implicit function theorem. Furthermore, either all 3 of them stay real-valued for $n\ne n^\ast$ locally or two of them split into a complex-conjugated non-real pair. It remains to prove that only the latter can be true: assume contrarily that $\kappa_{n,1},\kappa_{n,2},\kappa_{n,3}$ are real-valued trajectories for $n\in (n^\ast-\epsilon,n^\ast+\epsilon)$. By Rolle's theorem in case of $\kappa_{n,i}\ne \kappa_{n,j}$ for some $1\leq i,j\leq 3$ (let then $\kappa_{n,i}< \kappa_{n,j}$ without loss of generality) and by the product rule otherwise, there must be critical points $\kappa_{n,i,j}$ of the real-analytic function $\kappa \mapsto F_p(n,\kappa)$ for $n\in (n^\ast-\epsilon,n^\ast+\epsilon)$ such that
	\begin{align}\label{contra}
		\kappa_{n,i}\leq\kappa_{n,i,j}\leq \kappa_{n,j}\ ,
	\end{align} 
	in particular $\kappa_{n^\ast,i,j}=\kappa^\ast$. Since the roots of $J_p$ are discrete, however, we conclude for $\epsilon$ sufficiently small that 
	\begin{align*}
		\partial_\kappa F_p(n,\kappa_{n,i,j})=0\quad \Leftrightarrow \quad \kappa_{n,i,j}=\kappa^\ast \quad \text{or} \quad \kappa_{n,i,j}=\frac{\sqrt{n^\ast}\kappa^\ast}{\sqrt{n}}
	\end{align*}
	for $n\in (n^\ast-\epsilon,n^\ast+\epsilon)$ according to \Cref{deriv}. In both of the latter cases $\kappa_{n,i,j}$ is continuously-differentiable with $|\kappa'_{n^\ast,i,j}|<\infty$. This contradicts \Cref{contra} as $n\to n^\ast$ since $\lim_{n\to n^\ast}\kappa'_{n,i}=-\infty$ for all $1\leq i\leq 3$ by \Cref{angle} below.
\end{proof}

\noindent Next, we show that complex-valued ITE trajectories are limited to approach and leave DELs in only one angular direction depending on the half space of the complex plane they emerge from and move into, respectively. Arising in crossing complex-conjugated pairs then, there is an ambiguity of trajectorial ingoing and outgoing directions at DELs.  

\begin{lemma}\label{angle}
	Let $\kappa_n$ be a continuous ITE trajectory of the unit disk $D$ for $n\in (n^\ast-\epsilon,n^\ast+\epsilon)$ with $\epsilon>0$ sufficiently small such that $1\not\in (n^\ast-\epsilon,n^\ast+\epsilon)$ whose eigenfunction pairs $(v_n,w_n)$ are given by \Cref{eigf} for some fixed $p\in \mathbb{N}_0$.
	If $J_p(\kappa_{n^\ast})=0$, then it holds that $\lim_{n\to n^\ast}|\kappa'_n|=\infty$ and
	\begin{align}\label{dkappa}
		\kappa'_n=-\frac{(n\kappa_n^2-p^2)}{2n(n-1)\kappa_n}-\frac{\kappa_n J'_p(\kappa_n)^2}{2n(n-1)J_p(\kappa_n)^2}
	\end{align}
	as long as $J_p(\kappa_n)\ne 0$.
	Further, 
	\begin{align*}
		\lim_{n\searrow n^\ast}\arg(\kappa'_n)\ ,\lim_{n\nearrow n^\ast}\arg(\kappa'_n)\in \{\pm \pi/3,\pi\} \ \mod\ 2\pi
	\end{align*} 
	for $n^\ast>1$ and 
	\begin{align*}
		\lim_{n\searrow n^\ast}\arg(\kappa'_n)\ ,
		\lim_{n\nearrow n^\ast}\arg(\kappa'_n)\in\{0,\pm 2\pi/3\}  \ \mod\ 2\pi
	\end{align*} 
	for $0<n^\ast<1$. In particular, ingoing/outgoing angles at $\kappa_{n^\ast}$ which are integer multiples of $\pi$ correspond to real-valued trajectories, respectively.
\end{lemma}
\begin{proof}
	By the implicit function theorem, $\kappa_n$ is differentiable with 
	\begin{align*}
		\kappa'_n=-\frac{\partial_n F_p(n,\kappa)_{|(n,\kappa_n)}}{\partial_\kappa F_p(n,\kappa)_{|(n,\kappa_n)}}
	\end{align*}
	as long as $\partial_\kappa F_p(n,\kappa)_{|(n,\kappa_n)}$
	does not vanish for $\kappa=\kappa_n$. For this to happen, $\kappa_n$ must be a DEL by \Cref{delk} since $J_p(\kappa_n)=0$ if and only if $J_p(\sqrt{n}\kappa_n)=0$ according to the ITP boundary data in \Cref{ITP}. 
	Using $F_p(n,\kappa_n)=0$, we obtain
	\begin{align}\label{deln}
		\partial_n F_p(n,\kappa)_{|(n,\kappa_n)}=\frac{n\kappa_n^2-p^2}{2n}J_p(\kappa_n)J_p(\sqrt{n}\kappa_n)+\frac{\kappa_n^2}{2\sqrt{n}} J'_p(\kappa_n)J'_p(\sqrt{n}\kappa_n)
	\end{align}
	so that $\kappa'_n$ becomes
	\begin{align*}
		\kappa'_n&=-\frac{n\kappa_n^2-p^2}{2n(n-1)\kappa_n}-\frac{\kappa_n J'_p(\kappa_n)J'_p(\sqrt{n}\kappa_n)}{2\sqrt{n}(n-1)J_p(\kappa)J_p(\sqrt{n}\kappa)} \\
		&=-\frac{n\kappa_n^2-p^2}{2n(n-1)\kappa_n}-\frac{\kappa_n J'_p(\kappa_n)^2}{2n(n-1)J_p(\kappa_n)^2}\ .
	\end{align*}
	In this form, we see that $\lim_{n\nearrow n^\ast}|\kappa'_n|=\infty$ since $J_p(\kappa_{n^\ast})=0$ while $J'_p(\kappa_{n^\ast})\ne 0$.
	
	In order to determine the trajectorial ingoing/outgoing directions as $n\to n^\ast$ but circumvent the tangential blow-up behavior, we note that $J_p(\kappa)=(\kappa-\kappa_{n^\ast})h(\kappa)$, where  $h(\kappa_{n^\ast})=J'_p(\kappa_{n^\ast})\ne 0$. We can therefore rewrite \Cref{dkappa} in a neighborhood of $\kappa_{n^\ast}$ as
	\begin{align*}
		\left[(\kappa_n-\kappa_{n^\ast})^3\right]'=3\kappa'_n(\kappa_n-\kappa_{n^\ast})^2=-\frac{3(n\kappa_n^2-p^2)}{2n(n-1)\kappa_n}(\kappa_n-\kappa_{n^\ast})^2-\frac{3\kappa_n J'_p(\kappa_n)^2}{2n(n-1)h(\kappa_n)^2}\ .
	\end{align*}
	Since $\kappa_n$ is continuous by assumption, $n\mapsto (\kappa_n-\kappa_{n^\ast})^3\in C^1((n^\ast-\epsilon,n^\ast+\epsilon))$. Note that
	\begin{align}\label{pm}
		[(\kappa_n-\kappa_{n^\ast})^3]'_{|n=n^\ast}=\lim_{n\to n\ast} -\frac{3\kappa_n J'_p(\kappa_n)^2}{2n(n-1)h(\kappa_n)^2}=-\frac{3\kappa_{n^\ast}}{2n^\ast(n^\ast-1)}\in \mathbb{R}\backslash \{0\}
	\end{align}
	whose sign depends on whether $n^\ast>1$ or $0<n^\ast<1$.
	By continuity of the complex argument modulo $2\pi$ in $\mathbb{C}\backslash \{0\}$ and the identity $\arg(z_1z_2)=\arg(z_1)\arg(z_2)$ modulo $2\pi$ for all $z_1,z_2\in \mathbb{C}\backslash \{0\}$, we obtain in case of $n^\ast>1$
	\begin{align}\label{v1}
		\begin{split}
			\lim_{n\to n^\ast}3\arg\big((\kappa_n-\kappa_{n^\ast})/(n-n^\ast)\big)&=\lim_{n\to n^\ast}\arg\big((\kappa_n-\kappa_{n^\ast})^3/(n-n^\ast)\big)\\ &=\arg\big([(\kappa_n-\kappa_{n^\ast})^3]'_{|n=n^\ast}\big)\\
			&=-\pi \ \mod\ 2\pi\ .
		\end{split}
	\end{align}
	Using $[(\kappa_n-\kappa_{n^\ast})^3]'=3(\kappa_n-\kappa_{n^\ast})^2\kappa'_n$ for $n\ne n^\ast$, we also have 
	\begin{align}\label{v2}
		\begin{split}
			\lim_{n\to n^\ast}\left[2\arg\big((\kappa_n-\kappa_{n^\ast})/(n-n^\ast)\big)+\arg(\kappa'_n)\right]&=\lim_{n\to n^\ast}\arg\big((\kappa_n-\kappa_{n^\ast})^2/(n-n^\ast)^2\kappa'_n\big)\\ 
			&=\lim_{n\to n^\ast}\arg\big(3(\kappa_n-\kappa_{n^\ast})^2\kappa'_n\big)\\
			&=\arg\big([(\kappa_n-\kappa_{n^\ast})^3]'_{|n=n^\ast}\big)\\
			&=-\pi \ \mod\ 2\pi
		\end{split}
	\end{align}
	for $n^\ast>1$. By continuity of $\kappa_n$, which then enforces convergence towards $\kappa_{n^\ast}$ along only one angular directions as $n\nearrow n^\ast$ or $n\searrow n^\ast$, we conclude that
	\begin{align*}
		\lim_{n\nearrow n^\ast}\arg(\kappa'_n)=\lim_{n\nearrow n^\ast}\arg\big((\kappa_n-\kappa_{n^\ast})/(n-n^\ast)\big)\in \{\pm \pi/3,\pi\} \ \mod\ 2\pi
	\end{align*}
	and 
	\begin{align*}
		\lim_{n\searrow n^\ast}\arg(\kappa'_n)=\lim_{n\searrow n^\ast}\arg\big((\kappa_n-\kappa_{n^\ast})/(n-n^\ast)\big)\in \{\pm \pi/3,\pi\} \ \mod\ 2\pi
	\end{align*}
	for $n^\ast>1$. Here, the admissible limit set in the two equations above is a consequence of \Cref{v1} and the preceding equality follows in comparison with \Cref{v2} each. For $0<n^\ast<1$, we proceed similarly, using that \Cref{pm} is negative in this case, and obtain 
	\begin{align*}
		\lim_{n\nearrow n^\ast}\arg(\kappa'_n)=\lim_{n\nearrow n^\ast}\arg\big((\kappa_n-\kappa_{n^\ast})/(n-n^\ast)\big)\in \{0,\pm 2\pi/3\} \ \mod\ 2\pi
	\end{align*}
	as well as 
	\begin{align*}
		\lim_{n\searrow n^\ast}\arg(\kappa'_n)=\lim_{n\searrow n^\ast}\arg\big((\kappa_n-\kappa_{n^\ast})/(n-n^\ast)\big)\in \{0,\pm 2\pi/3\} \ \mod\ 2\pi\ .
	\end{align*}		
	It remains to prove that the two angles which are integer multiples of $\pi$ correspond to locally real-valued trajectories, respectively. This follows, restricting to $n\nearrow n^\ast$ for the sake of presentation, if we can show that $\Im(\kappa'_n)\gtrless 0$ for all $\kappa_n$ with $\Im(\kappa_n)\gtrless 0$ sufficiently close to $\kappa_{n^\ast}$ such that $\Re(\kappa_n-\kappa_{n^\ast})>\vert \Im(\kappa_n)\vert$ in case of $n^\ast>n>1$ and $\Re(\kappa_n-\kappa_{n^\ast})<-\vert \Im(\kappa_n)\vert$ for $0<n<n^\ast<1$, cf. \Cref{flow} for better illustration. 
	Integration with respect to $n$ then implies $\lim_{n\nearrow n^\ast}\Im(\kappa_n)\gtrless 0$ and thus contradicts $\lim_{n\to n^\ast}\kappa_n \ne \kappa_{n^\ast}\in \mathbb{R}$.
	The complementing case $n\searrow n^\ast$ can be proven along the same lines, showing that ITE trajectories $\kappa_n$ cannot escape from the real axis tangentially into the complex plane.
	Our proofs are inspired by the fact that poles of order 2 induce locally hyperbolic sectors for autonomous holomorphic flows, see \cite{broughan2003holomorphic}.
	\begin{figure}[htbp]
		\centering
		\includegraphics[width=12cm]{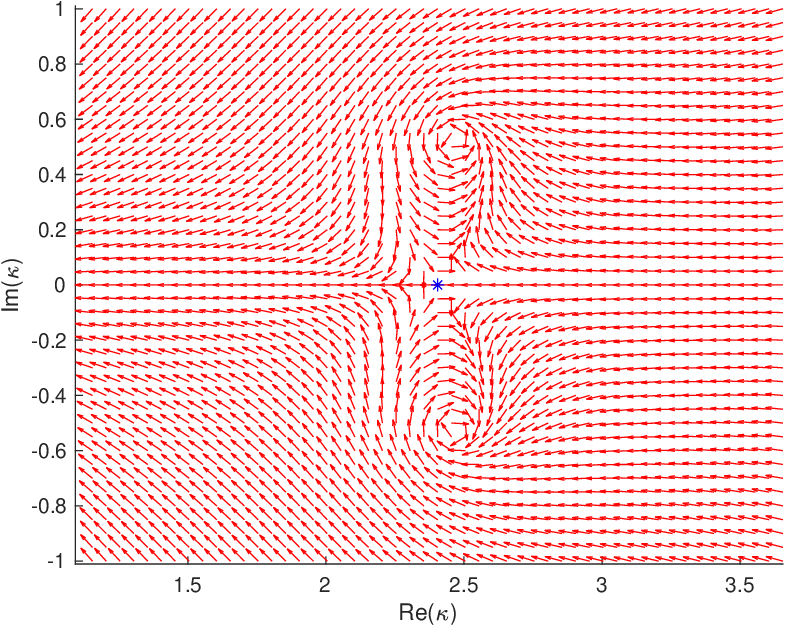}
		\caption{\label{flow}Snapshot of the holomorphic flow (red arrows) generated by $\kappa\mapsto d_n(\kappa)/|d_n(\kappa)|$ for $n=4$ and $p=0$. The plot is to illustrate that corresponding solutions, including ITE trajectories, cannot approach or escape the real axis tangentially from or into the complex plane near DELs (blue asterisk), respectively.} 
	\end{figure}

	In order to verify that $\Im(\kappa'_n)\gtrless 0$ for $n<n^\ast$ sufficiently close and all $\Im(\kappa_n)\gtrless 0$ subject to the aforementioned constraints, we confine ourselves to $\Im(\kappa_n)>0$ since ITEs occur in complex-conjugated pairs. We define for fixed $n$ the real-meromorphic function
	\begin{align*}
		d_n(\kappa):=-\frac{n\kappa^2-p^2}{2n(n-1)\kappa}-\frac{g_p(\kappa)}{2n(n-1)(\kappa-\kappa_{n^\ast})^2}
	\end{align*}
	in terms of \Cref{dkappa}, where
	\begin{align*}
		g_p(\kappa):=\frac{\kappa J'_p(\kappa)^2}{h(\kappa)^2}\ .
	\end{align*}
	We show in the following that $\arg(d_n(\kappa))$ is dominated by 
	\begin{align}\label{dom}
		\begin{split}
			&\arg\big(-(n-1)^{-1}(\kappa-\kappa_{n^\ast})^{-2}\big)\\
			=&
			\bigg\{
			\begin{array}{ll}
				\pi-2\arg(\kappa-\kappa_{n^\ast})\ , & \text{if }n>1,\ \Re(\kappa-\kappa_{n^\ast})> \Im(\kappa)>0\ ,\\
				-2\arg(\kappa-\kappa_{n^\ast})\ , & \text{if }0<n<1,\ \Re(\kappa-\kappa_{n^\ast})< -\Im(\kappa)<0
			\end{array}
			\ \mod\ 2\pi
		\end{split}
	\end{align}
	in a neighborhood of $\kappa_{n^\ast}$ which then yields $\Im(\kappa_n)>0$ since $\Im(-(n-1)^{-1}(\kappa-\kappa_{n^\ast})^{-2})> 0$ for all $\kappa$ under consideration. In order to make the argument rigorous, we observe that $g_p(\kappa_{n^\ast})=\kappa_{n^\ast}>0$, and a straighforward yet lengthy calculation shows that $g_p'(\kappa_{n^\ast})=0$ as well as 
	\begin{align}\label{gp2nd}
		g_p''(\kappa_{n^\ast})=\frac{-8\kappa_{n^\ast}^2+8p^2+1}{6\kappa_{n^\ast}}<0\ .
	\end{align}
	The last inequality follows by standard lower bounds on the roots of Bessel functions, see for instance \cite[p. 486]{Wa66}.	The second-order approximation of $g_p$ around $\kappa_{n^\ast}$ therefore reveals for $\kappa$ sufficiently close to $\kappa_{n^\ast}$ that
	\begin{align*}
		\Im(g_p(\kappa))=
		\left\{
		\begin{array}{ll}
			\leq 0\ , & \text{if }\Re(\kappa-\kappa_{n^\ast})> \Im(\kappa)>0\ ,\\
			\geq 0\ , & \text{if }\Re(\kappa-\kappa_{n^\ast})<- \Im(\kappa)<0\ .
		\end{array}
		\right.
	\end{align*}
	In the first case of the right-hand side above we can estimate for $n^\ast>n>1$
	\begin{align}\label{b1}
		\begin{split}
			0<\frac{\pi}{2}-\theta &<\arg\left(-\frac{g_p(\kappa)}{2n(n-1)(\kappa-\kappa_{n^\ast})^2}\right)\\
			&< \arg\big(-(n-1)^{-1}(\kappa-\kappa_{n^\ast})^{-2}\big)<\pi\  \ \mod\ 2\pi,
		\end{split}
	\end{align}
	where we restrict additionally $\kappa$ such that
	$-\theta<\arg(g_p(\kappa))\leq 0$ modulo $2\pi$
	for any fixed $0<\theta<\pi/2$. Note that $\theta$ can indeed be chosen arbitrarily small in a neighborhood around $\kappa_{n^\ast}$ since $g_p(\kappa_{n^\ast})>0$. In the other case, that is  $\Re(\kappa-\kappa_{n^\ast})<- \Im(\kappa)<0$ and $\kappa$ sufficiently close to $\kappa_{n^\ast}$, we obtain similarly for $0<n<n^\ast<1$
	\begin{align*}
		0&<\arg\big(-(n-1)^{-1}(\kappa-\kappa_{n^\ast})^{-2}\big)\\
		&<\arg\left(-\frac{g_p(\kappa)}{2n(n-1)(\kappa-\kappa_{n^\ast})^2}\right)<\frac{\pi}{2}+\theta <\pi\ \mod\ 2\pi\ .
	\end{align*}
	In both cases, we see that 
	\begin{align}\label{lhs}
		\Im\left(-\frac{g_p(\kappa)}{2n(n-1)(\kappa-\kappa_{n^\ast})^2}\right)>0
	\end{align}
	and, to verify that $\Im(d_n(\kappa))$ is likewise positive, we are left to show that adding $-(n\kappa^2-p^2)/(2n(n-1)\kappa)$ to the left-hand side of \Cref{lhs} does not affect the sign. 
	We first compute for $1<n<n^\ast$ and $\Re(\kappa-\kappa_{n^\ast})> \Im(\kappa)>0$ 
	\begin{align}\label{a1}
		\begin{split}
			\pi
			&<\pi+\arg\left(\frac{n\kappa^2-p^2}{2n(n-1)\kappa}\right)<\pi+\arg(n\kappa^2-p^2)<\pi+\arg(\kappa^2-p^2)\\
			&<\pi+\arg(\kappa^2-\kappa_{n^\ast}^2)=\pi+\arg\big((\kappa-\kappa_{n^\ast})(\kappa+\kappa_{n^\ast})\big)\\
			&<\pi+2\arg(\kappa-\kappa_{n^\ast})<\frac{3\pi}{2}\ \ \mod\ 2\pi\ .
		\end{split}
	\end{align}
	Here, the first inequality from the bottom line holds since $p<\kappa_{n^\ast}$ as a root of $J_p$, see \cite[p. 485]{Wa66}. Combining \Cref{b1} and \Cref{a1}, we obtain for $\arg(\kappa-\kappa_{n^\ast})>0$ modulo $2\pi$ small that
	\begin{align*}
		\Im(d_n(\kappa))>&\left\vert\frac{n\kappa^2-p^2}{2n(n-1)\kappa}\right\vert\Im\left(\mathrm{e}^{\mathrm{i}(\pi+2\arg(\kappa-\kappa_{n^\ast}))}\right)\\ &+
		\left\vert\frac{g_p(\kappa)}{2n(n-1)(\kappa-\kappa_{n^\ast})^2}\right\vert
		\Im\left(\mathrm{e}^{\mathrm{i}(\pi-2\arg(\kappa-\kappa_{n^\ast}))}\right)\ .	
	\end{align*}
	Since the first summand is bounded as $\kappa\to \kappa_{n^\ast}$ while the second blows up, we can find a neighborhood of $\kappa_{n^\ast}$ such that 
	\begin{align*}
		\Im(d_n(\kappa))>\left(\left\vert\frac{g_p(\kappa)}{2n(n-1)(\kappa-\kappa_{n^\ast})^2}\right\vert-\left\vert\frac{n\kappa^2-p^2}{2n(n-1)\kappa}\right\vert\right)
		\Im\left(\mathrm{e}^{\mathrm{i}(\pi-2\arg(\kappa-\kappa_{n^\ast}))}\right)>0 	
	\end{align*}
	which thus yields the assertion in the case $n^\ast>n>1$ and $\kappa$ sufficiently close to $\kappa_{n^\ast}$ such that $\Re(\kappa-\kappa_{n^\ast})> \Im(\kappa)>0$. In the other case, that is $\Re(\kappa-\kappa_{n^\ast})<- \Im(\kappa)<0$ and $0<n<n^\ast<1$, we even have 
	\begin{align}\label{easy}
		\Im\left(-\frac{n\kappa^2-p^2}{2n(n-1)\kappa}\right)>0
	\end{align}
	as long as $\Re(n\kappa^2)>p^2$ which is appropriate to assume since $\sqrt{n^\ast}\kappa_{n^\ast}>p$ as a DEL according to \Cref{cor} and \cite[p. 485]{Wa66}. In particular, \Cref{lhs} ensures that also $\Im(d_n(\kappa))>0$ for $0<n<n^\ast<1$ and $\kappa$ sufficiently close to $\kappa_{n^\ast}$ such that $\Re(\kappa-\kappa_{n^\ast})< -\Im(\kappa)<0$. This finally completes the proof of $\lim_{n\nearrow n^\ast}\arg(\kappa_n')\not \in \{0,\pi\}$ modulo $2\pi$ for complex-valued ITE trajectories $\kappa_n\to \kappa_{n^\ast}$. 
\end{proof}


\noindent We can also specify the geometric multiplicity of ITEs of the unit disk $\kappa_{n^\ast}$ for $n^\ast\ne 1$ such that $J_p(\kappa_{n^\ast})=0$. Depending on the Bessel index $p$, it is at least 3 for $p=0$, 5 for $p=1$ and 6 for $p>1$. Note a multiplicity of 1 for $p=0$ and 2 for integer $p>1$ would have been a priori expected according to \Cref{eigf}.
The increase originates from additional ITE trajectories $\widetilde{\kappa}_n$ with Bessel order $\widetilde{p}\in \{p-1,p+1\}$ which also pass the DEL as $n\to n^\ast$ along the real axis. The real-axis-constraint follows by \Cref{difference} and $J_{\widetilde{p}}(\kappa_{n^\ast})\ne 0$, cf. Bourget's hypothesis in \cite[p. 484]{Wa66}.

\begin{lemma}\label{multiplicity}
	Let $\kappa_{n^\ast}$ be an ITE of the unit disk $D$ for some $1\ne n^\ast>0$ such that $J_p(\kappa_{n^\ast})=0$. Then 
	\begin{align*}
		v_q(r,\phi)&=J_q(\kappa_{n^\ast} r)\cos(q\phi)\ , \quad \Big(v_q(r,\phi)=J_q(\kappa_{n^\ast} r)\sin(q\phi)\ , q>0\Big)\ , \\
		w_q(r,\phi)&=\alpha_n\sqrt{n^\ast} J_q(\sqrt{n^\ast}\kappa_{n^\ast} r)\cos(q\phi)\ , \quad \Big(w_q(r,\phi)=\alpha_n\sqrt{n^\ast}J_q(\sqrt{n}\kappa_{n^\ast} r)\sin(q\phi)\ , q> 0\Big),
	\end{align*}
	with $q\in \{p-1,p+1\}$ for $p>0$ and $q=1$ for $p=0$, are further linear independent ITP eigenfunction pairs for the ITE $\kappa_{n^\ast}$ besides \Cref{eigf}.
\end{lemma}
\begin{proof}
	By \Cref{big} there exists a real-valued continuous ITE trajectory $\kappa_n$ in a neighborhood of $n^\ast$ which is differentiable for $n\ne n^\ast$.
	Consider the limit quotients 
	\begin{align*}
		\widetilde{v}_q(r,\phi):=&\lim_{n\nearrow n^\ast}\frac{J_p(\kappa_{n}r)-J_p(\kappa_{n^\ast}r)}{\kappa_n-\kappa_{n^\ast}}\cos(q\phi)=J'_p(r\kappa_{n^\ast})r\cos(q\phi)\ , \\
		\widetilde{w}_q(r,\phi):=&\lim_{n\nearrow n\ast}\frac{\alpha_nJ_p(\sqrt{n}\kappa_{n}r)-\alpha_{n^\ast}J_p(\sqrt{n^\ast}\kappa_{n^\ast}r)}{\kappa_n-\kappa_{n^\ast}}\cos(q\phi)\\
		=&\ \alpha_{n^\ast}J'_p(r\sqrt{n^\ast}\kappa_{n^\ast})r\sqrt{n^\ast}\cos(q\phi)\\ 
		&\quad +\lim_{n\nearrow n\ast}\left(\alpha_{n}J'_p(r\sqrt{n}\kappa_{n})r\frac{\kappa_{n}}{2\sqrt{n}\kappa'_{n}}+\frac{\alpha'_{n}}{\kappa'_{n}}J_p(r\sqrt{n}\kappa_{n})\right)\cos(q\phi)\ ,
	\end{align*}
	whose right-hand sides are obtained by l'Hospital's rule. Here, $\alpha_n'$	denotes the derivative of
	\begin{align}\label{alpha}
		\alpha_n=\frac{J'_p(\kappa_n)}{\sqrt{n}J'_p(\sqrt{n}\kappa_n)}\ .
	\end{align}
	This representation is indeed well-defined in a neighborhood of $n^\ast$ since the roots of $J_p$ and $J'_p$ are distinct and $J_p(\kappa_{n^\ast})=J_p(\sqrt{n^\ast}\kappa_{n^\ast})=0$ by assumption.
	We thus obtain
	\begin{align*}
		\alpha'_{n}&=\frac{J''_p(\kappa_{n})\kappa'_{n}\sqrt{n}J'_p(\sqrt{n}\kappa_{n})-J'_p(\kappa_{n})\left[\frac{1}{2\sqrt{n}}J'_p(\sqrt{n}\kappa_{n})+\sqrt{n}J''_p(\sqrt{n}\kappa_{n})\left(\frac{\kappa_{n}}{2\sqrt{n}}+\sqrt{n}\kappa'_{n}\right)\right]}{nJ'_p(\sqrt{n}\kappa_{n})^2}	\\ 
		&=(n-1)J'_p(\kappa_n)J_p(\sqrt{n}\kappa_n)\kappa'_n\ ,
	\end{align*}
	where we replaced the second-order derivatives of $J_p$ above via \Cref{bessel}.
	It follows from the blow-up of $\lim_{n\to n^\ast}|\kappa'_n|=\infty$, cf. \Cref{angle}, that 
	\begin{align*}
		\lim_{n\nearrow n\ast}\alpha_{n}J'_p(r\sqrt{n}\kappa_{n})r\frac{\kappa_{n}}{2\sqrt{n}\kappa'_{n}}+\frac{\alpha'_{n}}{\kappa'_{n}}J_p(r\sqrt{n}\kappa_{n})=0
	\end{align*}
	within the definition of $\widetilde{w}_q(r,\phi)$.
	By construction, we have that both the Dirichlet data of $\widetilde{v}_q$ and $\widetilde{w}_q$ coincide as well as their Neumann data.
	The same also holds for 
	\begin{align*}
		\widehat{v}_q(r,\phi)&:=\frac{2p}{r\kappa_{n^\ast}}J_p(r\kappa_{n^\ast})\cos(q\phi)\ ,\\
		\widehat{w}_q(r,\phi)&:=\alpha_{n^\ast}\sqrt{n^\ast}\frac{2p}{r\sqrt{n^\ast}\kappa_{n^\ast}}J_p(r\sqrt{n^\ast}\kappa_{n^\ast})\cos(q\phi)\ .
	\end{align*}	
	However, neither $(\widetilde{v}_q,\widetilde{w}_q)$ nor $(\widehat{v}_q,\widehat{w}_q)$ satisfy the PDE conditions of \Cref{ITP}. Using the recurrence relations
	\begin{align*}
		2J'_p(z)=J_{p-1}(z)-J_{p+1}(z) \quad \text{and} \quad  \frac{2p}{z}J_p(z)=J_{p-1}(z)+J_{p+1}(z)\ ,
	\end{align*}
	we find that both
	\begin{align*}
		v_{p+1}(r,\phi)&:=\frac{\widehat{v}_{p+1}(r,\phi)-\widetilde{v}_{p+1}(r,\phi)}{2}=J_{p+1}(\kappa_{n^\ast}r)\cos((p+1)\phi)\ ,\\
		w_{p+1}(r,\phi)&:=\frac{\widehat{w}_{p+1}(r,\phi)-\widetilde{w}_{p+1}(r,\phi)}{2}=\alpha_{n^\ast}\sqrt{n^\ast}J_{p+1}(\sqrt{n^\ast}\kappa_{n^\ast}r)\cos((p+1)\phi)
	\end{align*}
	and
	\begin{align*}
		v_{p-1}(r,\phi)&:=\frac{\widehat{v}_{p-1}(r,\phi)+\widetilde{v}_{p-1}(r,\phi)}{2}=J_{p-1}(\kappa_{n^\ast}r)\cos((p-1)\phi)\ ,\\
		w_{p-1}(r,\phi)&:=\frac{\widehat{w}_{p-1}(r,\phi)+\widetilde{w}_{p-1}(r,\phi)}{2}=\alpha_{n^\ast}\sqrt{n^\ast}J_{p-1}(\sqrt{n^\ast}\kappa_{n^\ast}r)\cos((p-1)\phi)
	\end{align*}
	are linear independent ITP eigenfunction pairs for the ITE $\kappa_{n^\ast}$ as long as $p\ne 0$. If $p=0$, we effectively generate only one additional ITE eigenfunction pair in this way since $J_{-1}=-J_{1}$. We can repeat the proof with sine in place of cosine as angular part, provided $q\ne 0$ which only affects the case $p=1$. 
\end{proof}

\noindent For the rest of this section, we turn our attention to the behavior of complex-valued ITE trajectories for large $n$. We start with a result for general scatterers $D$ on possible accumulation points of non-real ITE sequences as $n\to\infty$.

\begin{theorem}\label{big2}
	Let $D$ be a bounded domain with $C^{1,1}$-boundary and let $\kappa_{n_m}$ be a sequence of non-real ITEs of $D$ with $n_m\nearrow \infty$ as $m\to\infty$. Denote by $(v_{n_m},w_{n_m})$ corresponding ITP eigenfunction pairs which are assumed to be normalized such that $\|v_{n_m}\|_{L^2(D)}=1$ for all $m\in\mathbb{N}$. If $\lim_{m\to\infty}\kappa_{n_m}=\kappa_\infty$, then $\kappa_\infty$ is a DEL of $D$ and a subsequence of $\{v_{n_m}\}_{m\in\mathbb{N}}$ converges strongly to some eigenfunction of $\kappa_\infty$ in $L^2(D)$.
\end{theorem}
\begin{proof}
	The normalization $\|v_{n_m}\|_{L^2(D)}=1$ implies with  \Cref{zerolemma} that
	\begin{align}\label{bound}
		\|v_{n_m}\|^2_{L^2(D)}=n_m\|w_{n_m}\|^2_{L^2(D)}=1 
	\end{align} 
	for all $m\in\mathbb{N}$ and thus $\lim_{m\to\infty}\|w_{n_m}\|_{L^2(D)}=0$.
	To show that $u_{n_m}:=(v_{n_m}-w_{n_m})\in H^2_0(D)$ is even uniformly bounded in $H^1_0(D)$ with respect to $m$, we first compute 
	\begin{align*}
		&\|\nabla u_{n_m}\|^2_{L^2(D)}\\
		=&\ \int_{D}(-\Delta u_{n_m})\overline{u}_{n_m}\, \mathrm{d}x\\
		=&\ \int_{D}\kappa_n^2(v_{n_m}-n_mw_{n_m})(\overline{v}_{n_m}-\overline{w}_{n_m})\, \mathrm{d}x \\
		=&\ \kappa_{n_m}^2\int_{D}|v_{n_m}|^2+n_m|w_{n_m}|^2 - n_m\overline{v}_{n_m}w_{n_m}-v_{n_m}\overline{w}_{n_m}\, \mathrm{d}x \\
		=&\ \kappa_{n_m}^2\int_{D}|v_{n_m}|^2+n_m|w_{n_m}|^2 -2 v_{n_m}\overline{w}_{n_m}\, \mathrm{d}x + \kappa_{n_m}^2\int_{D} v_{n_m}\overline{w}_{n_m} - n_m\overline{v}_{n_m}w_{n_m}\, \mathrm{d}x \ .
	\end{align*}
	The magnitude of the first integral in the bottom line is uniformly bounded because of $\kappa_{n_m}\to \kappa_\infty$ and \Cref{bound}. The second integral has a critical dependence on $n_m$ at first glace, but it is purely imaginary according to
	\begin{align*}
		&\ \kappa_{n_m}^2\int_{D} v_{n_m}\overline{w}_{n_m} - n_m\overline{v}_{n_m}w_{n_m}\, \mathrm{d}x \\
		=&\int_{D} -\Delta v_{n_m}\overline{w}_{n_m} +  \overline{v}_{n_m}\Delta w_{n_m}\, \mathrm{d}x \\
		=&\int_{D} -\Delta v_{n_m}\overline{w}_{n_m} +  \overline{v}_{n_m}\Delta w_{n_m}\, \mathrm{d}x 
		+\int_{D} \Delta w_{n_m}\overline{w}_{n_m} -  \overline{w}_{n_m}\Delta w_{n_m}\, \mathrm{d}x\\
		=&\int_{D}\Delta w_{n_m} \overline{u}_{n_m} -  \overline{w}_{n_m}\Delta u_{n_m}\, \mathrm{d}x  \\
		=&\int_{D}w_{n_m} \Delta\overline{u}_{n_m} -  \Delta\overline{w}_{n_m} u_{n_m}\, \mathrm{d}x  \\
		=&\ \overline{\kappa}_{n_m}^2\int_{D} w_{n_m}(-\overline{v}_{n_m} + n_m \overline{w}_{n_m})-(-n_m\overline{w}_{n_m})(v_{n_m}-w_{n_m})\, \mathrm{d}x \\
		=&\ -\overline{\kappa_{n_m}^2\int_{D} v_{n_m}\overline{w}_{n_m} - n\overline{v}_{n_m}w_{n_m}\, \mathrm{d}x}\ .
	\end{align*}
	Note that the boundary terms from integration by parts vanish since $u_{n_m}\in H^2_0(D)$.
	Thus we obtain  
	\begin{align*}
		\limsup_{m\to\infty}\|\nabla u_{n_m}\|^2_{L^2(D)}=\limsup_{m\to\infty}\Re\left(\kappa_{n_m}^2\int_{D}|v_{n_m}|^2+n_m|w_{n_m}|^2 -2 v_{n_m}\overline{w}_{n_m}\, \mathrm{d}x\right)= 2\Re\big(\kappa_\infty^2\big)
	\end{align*}
	and infer that there exists a weakly convergent subsequence of $\{u_{n_m}\}_{m\in\mathbb{N}}$, which we do not relabel, such that $u_{n_m}\rightharpoonup u_\infty$ in $H^1_0(D)$ and $u_{n_m}\to u_\infty$ in $L^2(D)$ as compact embedding. Since $v_{n_m}=u_{n_m}+w_{n_m}$ and $w_{n_m}\to 0$ in $L^2(D)$, we also get that $v_{n_m}\to u_\infty$ in $L^2(D)$. In particular, $\|u_\infty\|_{L^2(D)}=1$ as a consequence of \Cref{bound} and $u_\infty$ is a distributional solution to the Helmholtz equation in $H^1_0(D)$ with wave number $k_\infty$. That is, $k_\infty\in \mathbb{R}\backslash \{0\}$ is a DEL with eigenfunction $u_\infty$.
\end{proof}

In order to apply the result from the previous theorem to complex-valued ITE trajectories, we need to transition from local to global domains of $\kappa_n$ with respect to $n$.
In case of the unit disk $D$ we can show that local ITE trajectories do admit global continuous representatives which exist for all larger $n<\infty$ and which have the same Bessel order, cf. \eqref{eigf}. The extension is even unique if one additionally imposes equality of ingoing and outgoing directions at trajectorial crossing points. With this convention, locally real-valued ITE trajecories stay real-valued as $n$ varies while locally complex-valued ITE trajectories stay complex-valued according to the following lemma.
\begin{lemma}\label{ext}
	Let $\kappa_n$ be a continuous ITE trajectory of the unit disk $D$ for $n\in (n_{\min},n_{\max})$ such that $1\leq n_{\min}<n_{\max}<\infty$ whose eigenfunction pairs $(v_n,w_n)$ are given by \Cref{eigf} for some fixed $p\in \mathbb{N}_0$. Then $\kappa_n$ can be extended to a continuous ITE trajectory on $(n_{\min},\infty)$ subject to the same $p$ which is bounded as $n\to\infty$. The extension is unique if we impose additionally that $n\mapsto\arg(\kappa_n')$ is continuous when intersecting zeros of $J_p$. 
\end{lemma}
\begin{proof}
	By the global version of the Picard-Lindelöf theorem it holds that if $\kappa_n$ cannot be extended over $n_{\max}<\infty$ as a continuously differentiable solution to \eqref{dkappa}, then 
	\begin{align}\label{blowup}
		\limsup_{n\nearrow n_{\max}}|\kappa_n|=\infty
	\end{align}	
	or 
	\begin{align}\label{singularity}
		\lim_{n\nearrow n_{\max}}\kappa_n=\kappa^\ast\ .
	\end{align}	
	Here, $\kappa^\ast$ is a root of $J_p$ and thus a singularity of the right-hand side of \eqref{dkappa}. 
	In case of \eqref{singularity} $\kappa_n$ still admits a continuous trajectorial extension across $n_{\max}$ which follows by \Cref{big} since $n_{\max}$ satisfies the form of \eqref{ratio} by \Cref{cor}. Hence, it remains to show, and to which the rest of the proof below will be devoted, that blow-ups like \eqref{blowup} do not occur for $n_{\max}<\infty$ which then yields the desired extended domain $(n_{\min},\infty)$ for each originally local ITE trajectory. The final statement on uniqueness of trajectorial extensions is then again a consequence of the Picard-Lindelöf theorem which guarantees uniqueness of solutions to \eqref{dkappa} as long as $\kappa_n$ is not a root of $J_p$, cf. \eqref{singularity}. When $\kappa_{n^\ast}$ becomes a root of $J_p$ for some $n^\ast<\infty$, \Cref{big} ensures a simultaneous intersection of three local ITE trajectories for that $n^\ast$ which can be paired to satisfy the continuity constraint for $n\mapsto\arg(\kappa_n')$ at $n^\ast$ each by \Cref{angle}, that is, $\lim_{n\to n^\ast}\arg(\kappa')=-\pi$ modulo $2\pi$ for real-valued ITE trajectories $\kappa_n$ and $\lim_{n\to n^\ast}\arg(\kappa')\in \{\pm \pi/3\}$ modulo $2\pi$ for complex-valued ITE trajectories. 
	
	In order to exclude the possibility of \eqref{blowup}, we prove that any $\kappa_n$ keeps bounded for growing $n$. For this purpose, we multiply both sides of \eqref{dkappa} with $\overline{\kappa_n}$ and take real parts to obtain 
	\begin{align}\label{ineq}
		\begin{split}
			\big[|\kappa_n|^2\big]'&=-\frac{n}{n(n-1)}|\kappa_n|^2 + \frac{p^2}{n(n-1)}\Re\left(\frac{\overline{\kappa_n}}{\kappa_n}\right)+\frac{1}{n(n-1)}\Re\left(-\frac{J_p'(\kappa_n)^2}{J_p(\kappa_n)^2}\right)|\kappa_n|^2\\
			&\leq -\frac{n}{n(n-1)}|\kappa_n|^2 + \frac{p^2}{n(n-1)}+\frac{1}{n(n-1)}\Re\left(-\frac{J_p'(\kappa_n)^2}{J_p(\kappa_n)^2}\right)|\kappa_n|^2\ .
		\end{split}
	\end{align}
	If we can show that the magnitude of the Bessel function quotient above is uniformly bounded by approximately $1$ along certain centered arcs in the complex plane with sufficiently large radii, then it holds that $\big[|\kappa_n|^2\big]'\leq 0$ for all $\kappa_n$ along  these arcs and all $n$ sufficiently large, making the first summand above the dominant term. In particular, once they lie inside, ITE trajectories cannot leave corresponding arcs any more as $n$ grows which thus proves the assertion. 
	
	To get better control over the Bessel function quotient in \eqref{ineq}, we recall from \cite[p. 199]{Wa66} the asymptotic expansion
	\begin{align}
		J_p(z)=\sqrt{\frac{2}{\pi z}}\left(\cos\left(z-\frac{p}{2}\pi-\frac{1}{4}\pi\right)+\mathrm{e}^{|\Im(z)|}\mathcal{O}\left(\frac{1}{|z|}\right)\right)
	\end{align}
	which is valid for large $|z|$ and $|\arg(z)|<\pi$ modulo $2\pi$. In the sequel, we assume that $|z|$ is large without explicit mention and we may even restrict to $0\leq\arg(z)\leq \pi/4$ modulo $2\pi$ according to the necessary ITE criterion $|\Re(\kappa_n)|>|\Im(\kappa_n)|$, see \cite{cakoni2010interior}, and since ITEs arise in complex-conjugated pairs. Combining the identities $2J'_p(z)=J_{p-1}(z)-J_{p+1}(z)$ and $\cos(z+\pi/2)-\cos(z-\pi/2)=2\cos(z+\pi/2)$ for all $z\in\mathbb{C}$, we obtain
	\begin{align*}
		\frac{J_p'(z)}{J_p(z)}=\frac{\cos\left(z-\frac{p-1}{2}\pi-\frac{1}{4}\pi\right)+\mathrm{e}^{\Im(z)}\mathcal{O}\left(\frac{1}{|z|}\right)}{\cos\left(z-\frac{p}{2}\pi-\frac{1}{4}\pi\right)+\mathrm{e}^{\Im(z)}\mathcal{O}\left(\frac{1}{|z|}\right)}\ .
	\end{align*}
	Inserting the definition of $\cos(z)=(\mathrm{e}^{-\mathrm{i}z}+\mathrm{e}^{\mathrm{i}z})/2$ and rearranging factors of exponential type, we get
	\begin{align}
		\begin{split}\label{Besselquot}
			\frac{J_p'(z)}{J_p(z)}&=\frac{\mathrm{e}^{-\mathrm{i}\left(\Re(z)-\frac{p-1}{2}\pi-\frac{1}{4}\pi\right)}\mathrm{e}^{\Im(z)}+\mathrm{e}^{\mathrm{i}\left(\Re(z)-\frac{p-1}{2}\pi-\frac{1}{4}\pi\right)}\mathrm{e}^{-\Im(z)}+\mathrm{e}^{\Im(z)}\mathcal{O}\left(\frac{1}{|z|}\right)}
			{\mathrm{e}^{-\mathrm{i}\left(\Re(z)-\frac{p}{2}\pi-\frac{1}{4}\pi\right)}\mathrm{e}^{\Im(z)}+\mathrm{e}^{\mathrm{i}\left(\Re(z)-\frac{p}{2}\pi-\frac{1}{4}\pi\right)}\mathrm{e}^{-\Im(z)}+\mathrm{e}^{-\Im(z)}\mathcal{O}\left(\frac{1}{|z|}\right)} \\
			&=\mathrm{e}^{\mathrm{i}\frac{\pi}{2}}\frac{1-\mathrm{e}^{2\mathrm{i}\left(\Re(z)-\frac{p}{2}\pi-\frac{1}{4}\pi\right)}\mathrm{e}^{-2\Im(z)}+\mathcal{O}\left(\frac{1}{|z|}\right)}
			{1+\mathrm{e}^{2\mathrm{i}\left(\Re(z)-\frac{p}{2}\pi-\frac{1}{4}\pi\right)}\mathrm{e}^{-2\Im(z)}+\mathcal{O}\left(\frac{1}{|z|}\right)}\ .
		\end{split}
	\end{align}
	Hence, we can conclude for now that
	\begin{align}\label{ubound}
		\left|\frac{J_p'(z)}{J_p(z)}\right|\leq \frac{\left|1-\cos\left(2\Re(z)-p\pi-\frac{1}{2}\pi\right)\mathrm{e}^{-2\Im(z)}\right|+\left|\sin\left(2\Re(z)-p\pi-\frac{1}{2}\pi\right)\mathrm{e}^{-2\Im(z)}\right|+\mathcal{O}\left(\frac{1}{|z|}\right)}
		{\left|1+\cos\left(2\Re(z)-p\pi-\frac{1}{2}\pi\right)\mathrm{e}^{-2\Im(z)}\right|-\left|\sin\left(2\Re(z)-p\pi-\frac{1}{2}\pi\right)\mathrm{e}^{-2\Im(z)}\right|-\mathcal{O}\left(\frac{1}{|z|}\right)}\ .
	\end{align}
	Next, we define concentric arc curves $z_s=z_s(t)$ for $t\in [0,\pi/4]$ by 
	\begin{align}\label{arc}
		z_s(t):=j_s' \mathrm{e}^{\mathrm{i}t}\ ,
	\end{align} 
	where $j_s':=(s+p/2-3/4)\pi$ and $s\in \mathbb{N}$ will be specified later. To bound \eqref{ubound} along $z_s$ we will treat the cases $t\in [0,\sqrt{\pi/(2j_s')}]$ and $t\in (\sqrt{\pi/(2j_s')},\pi/4]$ separately:
	
	\noindent For $t\in [0,\sqrt{\pi/(2j_s')}]$, we can bound the cosine argument in \eqref{ubound} with the help of the inequality $\cos(z)\geq 1-z^2/2$ for $z\in \mathbb{R}$ by
	\begin{align*}
		2j_s'-p\pi-\frac{1}{2}\pi\geq 2\Re(z_s(t))-p\pi-\frac{1}{2}\pi
		&=2j_s'\cos(t)-p\pi-\frac{1}{2}\pi \\
		&\geq 2j_s'\left(1-\frac{t^2}{2}\right)-p\pi-\frac{1}{2}\pi \\
		&\geq 2j_s'-p\pi-\frac{1}{2}\pi-\frac{\pi}{2} 
	\end{align*}
	so that
	\begin{align*}
		\cos\left(2\Re(z_s(t))-p\pi-\frac{1}{2}\pi\right)\geq 0
	\end{align*}
	implies
	\begin{align}\label{ubound2}
		\left|\frac{J_p'(z_s(t))}{J_p(z_s(t))}\right|\leq \frac{1+\left|\sin\left(2j_s'\cos(t)-p\pi-\frac{1}{2}\pi\right)\mathrm{e}^{-2j_s'\sin(t)}\right|+\mathcal{O}\left(\frac{1}{|j_s'|}\right)}
		{1-\left|\sin\left(2j_s'\cos(t)-p\pi-\frac{1}{2}\pi\right)\mathrm{e}^{-2j_s'\sin(t)}\right|-\mathcal{O}\left(\frac{1}{|j_s'|}\right)}
	\end{align}
	for $t\in [0,\sqrt{\pi/(2j_s')}]$. Similarly, using the bounds $t/2\leq \sin(t)\leq t$ and $\cos(t)\leq 1-t^2/4$ for corresponding $t$, we can estimate
	\begin{align*}
		\left|\sin\left(2j_s'\cos(t)-p\pi-\frac{1}{2}\pi\right)\mathrm{e}^{-2j_s'\sin(t)}\right|
		&\leq \left|\sin\left(2j_s'\cos(t)-p\pi-\frac{1}{2}\pi\right)\mathrm{e}^{-j_s't}\right| \\
		&\leq \left|\sin\left(2j_s'\left(1-\frac{t^2}{4}\right)-p\pi-\frac{1}{2}\pi\right)\mathrm{e}^{-j_s't}\right| \\
		&= \left|\sin\left(\frac{t^2}{2}\right)\mathrm{e}^{-j_s't}\right| \\
		&\leq  \frac{t^2}{2}\mathrm{e}^{-j_s't}\ .
	\end{align*}
	A simple calculation shows that the global maximum of $t\mapsto \big(\mathrm{e}^{-j_s't}t^2/2\big)$ for positive $t$ is given by $2\mathrm{e}^{-2}/j_s'$. Therefore, \eqref{ubound2} can be further simplified to  
	\begin{align}\label{res1}
		\left|\frac{J_p'(z_s(t))}{J_p(z_s(t))}\right|\leq \frac{1+\mathcal{O}\left(\frac{1}{|j_s'|}\right)}
		{1-\mathcal{O}\left(\frac{1}{|j_s'|}\right)}
	\end{align}
	for $t\in [0,\sqrt{\pi/(2j_s')}]$.
	
	\noindent In case of $t\in (\sqrt{\pi/(2j_s')},\pi/4]$, one easily verifies that 
	\begin{align*}
		\mathrm{e}^{-2j_s'\sin(t)}\leq \mathrm{e}^{-j_s't} \leq \mathrm{e}^{-j_s'\sqrt{\frac{\pi}{2j_s'}}}=\mathrm{e}^{-\sqrt{\frac{j_s'\pi}{2}}}
	\end{align*}
	which implies in virtue of \eqref{ubound}
	\begin{align}\label{res2}
		\left|\frac{J_p'(z_s(t))}{J_p(z_s(t))}\right|\leq	
		\frac{\left|1+\mathrm{e}^{-\sqrt{\frac{j_s'\pi}{2}}}\right|+\left|\mathrm{e}^{-\sqrt{\frac{j_s'\pi}{2}}}\right|+\mathcal{O}\left(\frac{1}{|j_s'|}\right)}
		{\left|1-\mathrm{e}^{-\sqrt{\frac{j_s'\pi}{2}}}\right|-\left|\mathrm{e}^{-\sqrt{\frac{j_s'\pi}{2}}}\right|-\mathcal{O}\left(\frac{1}{|j_s'|}\right)}\ .
	\end{align}
	Comparing \eqref{res1} and \eqref{res2}, we finally obtain
	\begin{align*}
		\left|\frac{J_p'(z_s(t))}{J_p(z_s(t))}\right|\leq \max\left\{\frac{1+\mathcal{O}\left(\frac{1}{|j_s'|}\right)}
		{1-\mathcal{O}\left(\frac{1}{|j_s'|}\right)},\frac{\left|1+\mathrm{e}^{-\sqrt{\frac{j_s'\pi}{2}}}\right|+\left|\mathrm{e}^{-\sqrt{\frac{j_s'\pi}{2}}}\right|+\mathcal{O}\left(\frac{1}{|j_s'|}\right)}
		{\left|1-\mathrm{e}^{-\sqrt{\frac{j_s'\pi}{2}}}\right|-\left|\mathrm{e}^{-\sqrt{\frac{j_s'\pi}{2}}}\right|-\mathcal{O}\left(\frac{1}{|j_s'|}\right)}\right\}
	\end{align*}
	for all $t\in [0,\pi/4]$, which can approximate 1 arbitrarily close by chossing $j_s'$, or equivalently $s$, sufficiently large in \eqref{arc}. 
\end{proof}

\noindent By the former lemma we can infer that any global complex-valued and continuous ITE trajectory of the unit disk converges to some DEL as $n\to\infty$. 

\begin{corollary}\label{converge}
	Let $\kappa_n$ be a continuous complex-valued ITE trajectory of the unit disk $D$ on $(n_{\min},\infty)$ for some $n_{\min}\geq 1$ whose eigenfunction pair $(v_n,w_n)$ is given by \Cref{eigf} for $p\in \mathbb{N}_0$ fixed. Then it holds that
	\begin{align*}
		\lim_{n\to\infty}\kappa_n=\kappa^\ast\ ,
	\end{align*}
	for some $\kappa^\ast$ such that $J_p(\kappa^\ast)=0$.
\end{corollary}
\begin{proof}
	By \Cref{ext} we know that $\kappa_n$ is bounded for large $n$. Hence, we can find a subsequence of $\kappa_n$ as $n\to\infty$ which converges to some DEL $\kappa^\ast$ of the unit disk according to \Cref{big2}. Also, a subsequence of corresponding ITP eigenfunctions $v_n$ converges strongly to some eigenfunction $v$ of $\kappa^\ast$. Since all $v_n$ are subject to the same Bessel order $p$ by assumption, $v$ itself must also have Bessel order $p$ which is due to orthogonality of Fourier Bessel functions in $L^2(D)$ for different integer indices $p$ as adopted by its angular Fourier basis. In particular, all accumulation points $\kappa^\ast$ of $\kappa_n$ are restricted to roots of $J_p$.
	
	\noindent Assume now contrarily that two subsequences of $\kappa_n$ converge to different roots $\kappa^\ast_1<\kappa^\ast_2$ of $J_p$ as $n\to\infty$. Choose any $\kappa_{1,2}$ which is not a root of $J_p$ such that $\kappa^\ast_1<\kappa_{1,2}<\kappa^\ast_2$. By assumption and continuity of $\kappa_n$, there is a sequence $n_m\to\infty$ as $m\to\infty$ such that $|\kappa_{n_m}|=\kappa_{1,2}$ and
	\begin{align}\label{contrad}
		[|\kappa_{n_m}|^2]'\geq 0
	\end{align}
	for all $m\in\mathbb{N}$.
	On the other hand, consider the sector 
	\begin{align*}
		z_{1,2}(t)=\kappa_{1,2}\mathrm{e}^{\mathrm{i}t}
	\end{align*}
	for $t\in [-\pi/4,\pi/4]$ (recall non-real ITEs $\kappa_n$ do not exist for $|\Re(\kappa_n)|<|\Im(\kappa_n)|$, see \cite{cakoni2010interior}) and set $C:=\max_{[-\pi/4,\pi/4]}|J_p'(z_{1,2}(t))/J_p(z_{1,2}(t))|^2<\infty$. 
	Since $p<\kappa^\ast_1<\kappa_{1,2}$, see \cite[p. 485]{Wa66}, we deduce $[|\kappa_{n_m}|^2]'< 0$ in \eqref{ineq} for all $n_m>C+1$ which is a contradiction to \eqref{contrad}.
\end{proof}

\begin{remark}
In contrast to complex-valued ITE trajectories, cf. \Cref{converge}, continuous and globally real-valued $\kappa_n$ fulfill $\lim_{n\to\infty}\kappa_n=0$. To see this, recall \eqref{ineq} which can be further estimated by
\begin{align*}
	[\kappa_n^2]'\leq -\frac{\kappa_n^2}{n-1}+\frac{p^2}{n(n-1)}\ ,
\end{align*}
where $p\in\mathbb{N}_0$ is again the underlying Bessel index according to \eqref{eigf}. Hence, $n\mapsto\kappa_n^2$ can only be non-decreasing if $\kappa_n^2\leq p^2/n$ which goes to zero as $n\to\infty$. Assuming contrarily that $\limsup_{n\to\infty}\kappa_n^2=\epsilon>0$, we infer that there is $n_\epsilon>1$ such that $n\mapsto \kappa_n^2$ is decreasing for all $n>n_\epsilon$, in particular $\lim_{n\to\infty}\kappa_n^2=\epsilon$ exists and $\kappa_n^2>\epsilon$ for all $n>n_\epsilon$. We obtain
\begin{align*}
	0&<\epsilon=\lim_{n\to\infty}\kappa_n^2\leq \kappa_{n_\epsilon}^2+\lim_{n\to\infty}\int_{n_\epsilon}^n[\kappa_t^2]'\,\mathrm{d}t\\
	&< \kappa_{n_\epsilon}^2+\limsup_{n\to\infty}\int_{n_\epsilon}^n-\frac{\epsilon}{t-1}+\frac{p^2}{t(t-1)}\,\mathrm{d}t<-\infty
\end{align*}
since $t\mapsto 1/(t-1)$ is not integrable in $(n_\epsilon,\infty)$, which yields a contradiction.
\end{remark}

\noindent Classifying continuous and complex-valued ITE trajectories by their unique global representative specified in \Cref{ext}, we can show that for any DEL there exists exactly one complex-conjugated pair of global complex-valued ITE trajectories which converge to that DEL as $n\to \infty$. In combination with \Cref{converge}, we may even conclude that there is a one-to-one correspondence between complex-valued ITE trajectories and DELs for the unit disk. This relation even holds including the geometric multiplicity of DELs and ITEs since they share the same eigenfunction structure, cf.  \eqref{eigf}, which results in a multiplicity of 1 for $p=0$ and 2 for $p>0$.

\begin{theorem}\label{theorem210}
	Let $\kappa^\ast$ be such that $J_p(\kappa^\ast)=0$. Then there exists one and only one complex-conjugated pair of global complex-valued and continuous ITE trajectories $(\kappa_n,\overline{\kappa}_n)$ on $(n_{\min},\infty)$ for some $n_{\min}\geq 1$ such that $n\mapsto\arg(\kappa_n')$ is continuous when intersecting zeros of $J_p$, whose eigenfunction pairs $(v_n,w_n)$ are given by \Cref{eigf} and which fulfills
	\begin{align*}
		\lim_{n\to\infty}\kappa_n=\kappa^\ast\ .
	\end{align*}
\end{theorem}
\begin{proof}
	We start with the existence part of the theorem. Let $\kappa^-<\kappa^\ast<\kappa^+$ be consecutive DELs of the unit disk or $\kappa^-=0$ in case that $\kappa^\ast$ is the smallest eigenvalue of the unit disk and fix any $a,b$ such that $\kappa^-<a<\kappa^\ast<b<\kappa^+$. Restricting to non-real ITEs with positive imaginary part which is feasible since ITEs arise in complex-conjugated pairs, we first want to show that for all $n$ sufficiently large there are no ITEs on both of the sectors 
	\begin{align}\label{sectors}
		z^-(t)=a\mathrm{e}^{\mathrm{i}t} \qquad \text{and} \qquad 
		z^+(t)=b\mathrm{e}^{\mathrm{i}t}
	\end{align}
	for $t\in(0,\pi/4]$ (recall that ITEs $\kappa_n$ do not exist for $\Re(\kappa_n)>\Im(\kappa_n)$, see \cite{cakoni2010interior}). To this end, note by \Cref{big2} that for any $\epsilon>0$ sufficiently small there exists $n_\epsilon>1$ such that there are no non-real ITEs $\kappa_n\in B_\epsilon(a)$ or $\kappa_n\in B_\epsilon(b)$ for all $n>n_\epsilon$, where $B_\epsilon(z)$ is the open ball around $z\in \mathbb{C}$ with radius $\epsilon$. Fix such an $\epsilon>0$ and set
	\begin{align*}
		C:=\max\left\{\max_{[\delta,\pi/4]}\left|\frac{J_p'(z^-(t))}{J_p(z^-(t))}\right|,\max_{[\delta,\pi/4]}\left|\frac{J_p'(z^+(t))}{J_p(z^+(t))}\right|\right\}<\infty\ ,
	\end{align*}
	where $\delta>0$ is chosen such that both $z^-(\delta)\in B_\epsilon(a)$ and $z^+(\delta)\in B_\epsilon(b)$. Using \eqref{Besselquot}, we can find $n_\delta>n_\epsilon$ such that 
	\begin{align*}
		\frac{C}{\sqrt{n}}<1-\delta<\min\left\{\min_{[\delta,\pi/4]}\left|\frac{J_p'(\sqrt{n}z^-(t))}{J_p(\sqrt{n}z^-(t))}\right|,\min_{[\delta,\pi/4]}\left|\frac{J_p'(\sqrt{n}z^+(t))}{J_p(\sqrt{n}z^+(t))}\right|\right\}
	\end{align*}
	for all $n>n_\delta$. Hence, the complementing assertion that the sectors $z^-(t)$ and $z^+(t)$ for all  $t\in [\delta,\pi/4]$ are also ITE-free if $n>n_\delta$ follows by \eqref{determinant} since
	\begin{align}\label{frac}
		F_p(n,\kappa_n)=0 \quad \Leftrightarrow \quad  \frac{J_p'(\kappa_n)}{\sqrt{n}J_p(\kappa_n)}=\frac{J_p'(\sqrt{n}\kappa_n)}{J_p(\sqrt{n}\kappa_n)}\ .
	\end{align}
	By \Cref{big} we can find some $n^\ast>n_\delta$ and a complex-conjugated pair of complex-valued continuous ITE trajectories $(\kappa_n,\overline{\kappa}_n)$ subject to Bessel order $p$ for $n\in (n^\ast-\epsilon,n^\ast+\epsilon)$ with $\epsilon>0$ such that $\kappa_{n^\ast}=\overline{\kappa}_{n^\ast}=\kappa^\ast$. We set $n_{\min}:=n^\ast-\epsilon$ and recall by \Cref{ext} that $\kappa_n$ has a global complex-valued and continuous extension for all $n>n_{\min}$ whose ingoing and outgoing directions coincide when intersecting the real axis. 
	By complex conjugation, the same holds true for $\overline{\kappa}_n$ which we thus do not treat separately any more in what follows. Finally, \Cref{converge} ensures that $\kappa_n$ converges to some DEL of the unit disk as $n\to\infty$. This limit must be $\kappa^\ast$ by construction of \eqref{sectors}.  
	
	\noindent Regarding the uniqueness assertion of our theorem, we will prove that there is $\gamma>0$ small such that for all $\widetilde{n}^\ast$ of the form \eqref{ratio} and sufficiently large there are no non-real ITEs in $B_\gamma(\kappa^\ast)$ with index of refraction $n=\widetilde{n}^\ast$ whose eigenfunction pairs have Bessel order $p$. Since any global complex-valued and continuous ITE trajectory $\widetilde{\kappa}_n$ that converges to $\kappa^\ast$ as $n\to\infty$ must lie inside of $B_\gamma(\kappa^\ast)$ for $n$ sufficiently large, but $\kappa^\ast$ is a root of order 3 of $\kappa\mapsto F_p(\widetilde{n}^\ast,\kappa)$ which allows for a single pair of intersecting complex-conjugated ITE trajectories at $n=\widetilde{n}^\ast$ only, cf. \Cref{big}, we conclude that $\widetilde{\kappa}_n=\kappa_n$ or $\widetilde{\kappa}_n=\overline{\kappa}_n$ for  $n\geq\widetilde{n}^\ast$ and thus also for $n>n_{\min}$. In particular, $\kappa_n=\widetilde{\kappa}_n=\kappa^\ast$ for infinitely many $n<\infty$.
	
	\begin{figure}[htbp]
		\centering
		\begin{tikzpicture}
			\draw [-latex,thick](1,0) -- (10,0);
			\draw [-latex,thick](1,-1) -- (1,3);
			\node[right, color=black] at (10, 0) (p1) { $\mathrm{Re}(\kappa)$ };
			\node[centered, color=black] at (4, -0.3) (p1) { ${\color{blue}{\kappa^\ast}}-\gamma$ };
			\node[centered, color=black] at (8, -0.3) (p1) { ${\color{blue}{\kappa^\ast}}+\gamma$ };
			\node[centered, color=blue] at (6, -0.3) (p1) { $\kappa^\ast$ };
			\node[centered, color=blue] at (6, -0.01) (p1) { $\ast$ };
			\node[centered, color=black] at (1, 3.3) (p1) { $\mathrm{Im}(\kappa)$ };
			
			\draw [-latex,thick,green](6.8750,1.5155) -- (6.8750+0.43301,1.5155-0.25);
			
			\draw[thick,dashed] (8,0) arc (0:60:2) ;
			\draw[thick,dashed] (6.3,0) arc (0:60:0.3) ;
			\draw[thick,red] ([shift=(60:2)]6,0) arc (60:120:2);
			\draw[thick,red] ([shift=(60:0.3)]6,0) arc (60:120:0.3);
			\draw[thick,dashed] ([shift=(120:2)]6,0) arc (120:180:2);
			\draw[thick,dashed] ([shift=(120:0.3)]6,0) arc (120:180:0.3);

			\draw[thick,green] ([shift=(240:0.2)]6.8750,1.5155) arc (240:330:0.2);
			
			\draw [thick,red](6.15,0.2598) -- (7,1.7321);
			\draw [thick,red](5.85,0.2598) -- (5,1.7321);
			
			\node[centered, color=red] at (6, 1.5) (p1) { $M_\gamma$ };
			\node[centered, color=green] at (4.75, 0.7) (p1) { $L_\gamma$ };
			\node[centered, color=green] at (7.25, 0.7) (p1) { $R_\gamma$ };
			\node[centered, color=green] at (7.8, 1.7) (p1) { $\mathrm{e}^{-\frac{\pi}{2}\mathrm{i}}v$ };
			\node[centered, color=green] at (6.9, 1.38) (p1) { $\cdotp$ };
			
			\draw [->,orange] plot [smooth] coordinates {(7.9,0.1) (7.4,0.13) (6.6,0.3) (6.72,1.1)   };
		\end{tikzpicture}
		\caption{\label{figX} Illustration of the sectors $L_\gamma$, $M_\gamma$, $R_\gamma$ which do not contain any non-real ITEs if the refractive index is such that the complex-conjugated pair of complex-valued ITE trajectories intersect the DEL $\kappa^\ast$.
			For $M_\gamma$ (red) this fact is shown by exploiting Rouch\'e's theorem in combination with \eqref{frac}, while for $L_\gamma$ or $R_\gamma$ (green) one uses that any complex-valued solution to \eqref{dkappa} (orange) would have hit the real axis or the dashed boundary line, both of which is excluded for ITE trajectories.}
	\end{figure}

	\noindent We will proceed in two steps to determine $\gamma>0$, cf. \Cref{figX} for better visualization of the following strategy: Restricting to the upper half plane of $\mathbb{C}$ again since ITEs arise in complex-conjugated pairs, we first show that the middle sector
	\begin{align*}
		M_{\gamma}:=\Big\{\kappa\in\mathbb{C}:\ \frac{\pi}{3}\leq\arg(\kappa-\kappa^\ast)\leq\frac{2\pi}{3}\quad \text{mod }2\pi, \  0<|\kappa-\kappa^\ast|<\gamma\Big\}
	\end{align*}
	does not contain ITEs for all $n=\widetilde{n}^\ast$ of the form \eqref{ratio} and sufficiently large if $\gamma>0$ is properly chosen independently of $n$. Then we investigate the outer sectors
	\begin{align*}
		R_{\gamma}&:=\Big\{\kappa\in\mathbb{C}:\ 0<\arg(\kappa-\kappa^\ast)<\frac{\pi}{3}\quad \text{mod }2\pi, \  0<|\kappa-\kappa^\ast|<\gamma\Big\}\ ,\\
		L_{\gamma}&:=\Big\{\kappa\in\mathbb{C}:\ \frac{2\pi}{3}<\arg(\kappa-\kappa^\ast)<\pi\quad \text{mod }2\pi, \  0<|\kappa-\kappa^\ast|<\gamma\Big\}\ ,
	\end{align*}
	correspondingly. We start with recalling Rouch\'e's theorem which states that $\kappa\mapsto F_p(\widetilde{n}^\ast,\kappa)$ does not have zeros in $M_{\gamma}\backslash B_\varepsilon(\kappa^\ast)$ for $\gamma>\varepsilon>0$ if
	\begin{align}\label{rouche}
		\left|\frac{J_p(\sqrt{\widetilde{n}^\ast}\kappa)}{\sqrt{\widetilde{n}^\ast}J_p'(\sqrt{\widetilde{n}^\ast}\kappa)}\right|<\left|\frac{J_p(\kappa)}{J_p'(\kappa)}\right| \qquad \text{for all } \kappa\in \partial [M_{\gamma}\backslash B_\varepsilon(\kappa^\ast)]\ .
	\end{align}
	By letting $\varepsilon\to 0$, we conclude with \eqref{frac} that $M_\gamma$ is ITE-free for $n=\widetilde{n}^\ast$ large and which are of the form \eqref{ratio}. In order to prove \eqref{rouche}, we define for $t\geq 0$
	\begin{align}\label{gn}
		g_n(t):=\frac{J_p(\sqrt{n}(\kappa^\ast+tv))}{\sqrt{n}J_p'(\sqrt{n}(\kappa^\ast+tv))}\ ,
	\end{align}
	where $v=\mathrm{e}^{\mathrm{i}\phi}$ and $\phi\in [\pi/3,2\pi/3]$ is an implicit parameter. Using \eqref{bessel} and $J_p(\kappa^\ast)=J_p(\sqrt{\widetilde{n}^\ast}\kappa^\ast)=0$ by definition of $\widetilde{n}^\ast$, straightforward calculations yield that
	\begin{align}\label{taylor1}
		\begin{split}
			g_1(0)&=g_{\widetilde{n}^\ast}(0)=0\ ,\\
			g_1'(0)&=g_{\widetilde{n}^\ast}'(0)=v\ ,\\
			g_1''(0)&=g_{\widetilde{n}^\ast}''(0)=\frac{v^2}{\kappa^\ast}\ .\\
		\end{split}
	\end{align}
	Taking derivatives on both sides of \eqref{bessel}, we further get
	\begin{align}\label{taylor2}
		g_1'''(0)=v^3\left(2-\frac{2p^2+3}{{\kappa^\ast}^2}\right),\qquad \text{and}\qquad g_{\widetilde{n}^\ast}'''(0)=\widetilde{n}^\ast v^3\left(2-\frac{2p^2+3}{\widetilde{n}^\ast{\kappa^\ast}^2}\right)
	\end{align}
	and by \cite[p. 486]{Wa66} we deduce that for any $\widetilde{n}^\ast>1$
	\begin{align*}
		0<g_1'''(0)/v^3<\widetilde{n}^\ast g_1'''(0)/v^3<g_{\widetilde{n}^\ast}'''(0)/v^3\ .
	\end{align*}
	We will consider $\partial [M_{\gamma}\backslash B_\varepsilon(\kappa^\ast)]$ piecewise and show first that there is $\widetilde{\varepsilon}>0$ (depending on $\widetilde{n}^\ast$) such that $|g_{\widetilde{n}^\ast}(\varepsilon)|<|g_1(\varepsilon)|$ for all $0<\varepsilon<\widetilde{\varepsilon}$ and $\phi\in [\pi/3,2\pi/3]$. We write
	\begin{align}\label{lower}
		\begin{split}
			|g_{\widetilde{n}^\ast}(\varepsilon)|^2&=|g_1(\varepsilon)+\big(g_{\widetilde{n}^\ast}(\varepsilon)-g_1(\varepsilon)\big)|^2\\
			&= |g_1(\varepsilon)|^2+2\Re\big(\overline{g_1(\varepsilon)}\big(g_{\widetilde{n}^\ast}(\varepsilon)-g_1(\varepsilon)\big)\big)+|g_{\widetilde{n}^\ast}(\varepsilon)-g_1(\varepsilon)|^2\ .
		\end{split}
	\end{align}
	and, since $g_n$ is a composition of holomorphic functions, we can assume that $\widetilde{\varepsilon}>0$ is such that both $|\arg(g_1(\varepsilon))-\arg(v)|<\pi/24$ modulo $2\pi$ with $|g_1(\varepsilon)|>\varepsilon/2$ and  $|\arg(g_{\widetilde{n}^\ast}'''(\varepsilon)-g_1'''(\varepsilon))-\arg(v^3)|<\pi/24$ modulo $2\pi$ with $2|g_{\widetilde{n}^\ast}'''(0)-g_1'''(0)|>|g_{\widetilde{n}^\ast}'''(\varepsilon)-g_1'''(\varepsilon)|>\Re\big((g_{\widetilde{n}^\ast}'''(\varepsilon)-g_1'''(\varepsilon))/v^3\big)>|g_{\widetilde{n}^\ast}'''(0)-g_1'''(0)|/2$ for all $0<\varepsilon<\widetilde{\varepsilon}$ and all $\phi\in [\pi/3,2\pi/3]$.
	By Taylor's theorem we have by \eqref{taylor1}
	\begin{align*}
		g_{\widetilde{n}^\ast}(\varepsilon)-g_1(\varepsilon)=\int_0^\varepsilon \big(g_{\widetilde{n}^\ast}'''(s)-g_1'''(s)\big)\frac{(\varepsilon-s)^2}{2}\,\mathrm{d}s\ ,
	\end{align*}
	hence for all $0<\varepsilon<\widetilde{\varepsilon}$
	\begin{align*}
		\frac{\varepsilon^3}{12}|g_{\widetilde{n}^\ast}'''(0)-g_1'''(0)|<\left|\frac{g_{\widetilde{n}^\ast}(\varepsilon)-g_1(\varepsilon)}{v^3}\right|=|g_{\widetilde{n}^\ast}(\varepsilon)-g_1(\varepsilon)|<\frac{\varepsilon^3}{3}|g_{\widetilde{n}^\ast}'''(0)-g_1'''(0)|\ .
	\end{align*}
	Due to $\arg(v^2)\in[2/3\pi,4/3\pi]$ modulo $2\pi$ by assumption on $v$, we also get that
	\begin{align*}
		&\Re\big(\overline{g_1(\varepsilon)}\big(g_{\widetilde{n}^\ast}(\varepsilon)-g_1(\varepsilon)\big)\big)\\
		=\,&|g_1(\varepsilon)|\cdotp|g_{\widetilde{n}^\ast}(\varepsilon)-g_1(\varepsilon)|\Re\left(\mathrm{e}^{\mathrm{i}(-\arg(g_1(\varepsilon))+\arg(g_{\widetilde{n}^\ast}(\varepsilon)-g_1(\varepsilon))}\right)\\
		=\,&|g_1(\varepsilon)|\cdotp|g_{\widetilde{n}^\ast}(\varepsilon)-g_1(\varepsilon)|\Re\left(\mathrm{e}^{\mathrm{i}(-\arg(v)+(\arg(v)-\arg(g_1(\varepsilon)))+\arg(v^3)+(\arg(g_{\widetilde{n}^\ast}(\varepsilon)-g_1(\varepsilon))-\arg(v^3))}\right)\\
		<\,&|g_1(\varepsilon)|\cdotp|g_{\widetilde{n}^\ast}(\varepsilon)-g_1(\varepsilon)|\Re\left(\mathrm{e}^{\mathrm{i}(2\phi-\pi/12)}\right)\\
		<\,&\frac{\varepsilon}{2}\cdotp\frac{\varepsilon^3}{12}\Re\left(\mathrm{e}^{\mathrm{i}\frac{7}{12}\pi}\right) = \frac{\varepsilon^4}{24}\Re\left(\mathrm{e}^{\mathrm{i}\frac{7}{12}\pi}\right)\ .
	\end{align*}
	Therefore, \eqref{lower} can be estimated for $0<\varepsilon<\widetilde{\varepsilon}$ by
	\begin{align*}
		|g_{\widetilde{n}^\ast}(\varepsilon)|^2< |g_1(\varepsilon)|^2+\frac{\varepsilon^4}{12}\Re\left(\mathrm{e}^{\mathrm{i}\frac{7}{12}\pi}\right)+\frac{\varepsilon^6}{9}|g_{\widetilde{n}^\ast}'''(0)-g_1'''(0)|^2 
	\end{align*}
	and since $\Re\left(\mathrm{e}^{\mathrm{i}\frac{7}{12}\pi}\right)<0$, we conclude by comparing powers of $\varepsilon$ that $\widetilde{\varepsilon}>0$ can additionally be chosen such that $|g_{\widetilde{n}^\ast}(\varepsilon)|<|g_1(\varepsilon)|$ for all $0<\varepsilon<\widetilde{\varepsilon}$ and all $\phi\in [\pi/3,2\pi/3]$.
	
	\noindent Considering the angular boundaries of $M_\gamma$ next, we need to show that $|g_{\widetilde{n}^\ast}(t)|<|g_1(t)|$ for all $0<t<\gamma$ with $\gamma>0$ to be defined, for all $\widetilde{n}^\ast$ sufficiently large and $\phi=\pi/3$ and $\phi=2\pi/3$ fixed, respectively. Since $g_{\widetilde{n}^\ast}$ solves the non-linear first order differential equation \eqref{ode} for $n=\widetilde{n}^\ast$ of the form \eqref{ratio} which can be proven by exploiting  \eqref{bessel}, existence of some $\gamma>0$ and a corresponding lower threshold $n_\gamma>1$ for $\widetilde{n}^\ast$ directly follows by \Cref{prop}. 
	
	\noindent To finish our proof of \eqref{rouche}, it remains to prove that $|g_{\widetilde{n}^\ast}(\gamma)|<|g_1(\gamma)|$ for all $\phi\in [\pi/3,2\pi/3]$ and $\widetilde{n}^\ast$ sufficiently large. To this end, recall from \eqref{frac} that $\kappa_{\widetilde{n}^\ast}=\gamma v$ is an ITE if and only if $g_{\widetilde{n}^\ast}(\gamma)=g_1(\gamma)$ for some $\phi\in [\pi/3,2\pi/3]$. We have that $\min_{\phi\in [\pi/3,2\pi/3]}|g_1(\gamma)|>0$ as the roots of $J_p$ are real-valued for $p\in\mathbb{N}_0$ and we know, comparing  \eqref{gn} with \eqref{Besselquot}, that $\lim_{\widetilde{n}^\ast\to\infty}|g_{\widetilde{n}^\ast}(\gamma)|=0$ uniformly in $\phi\in [\pi/3,2\pi/3]$ for all $\gamma>0$ fixed. Hence, given $\gamma>0$, it holds that $g_{\widetilde{n}^\ast}(\gamma)<g_1(\gamma)$ for all $\phi\in [\pi/3,2\pi/3]$ and $\widetilde{n}^\ast\geq n_\gamma$ sufficiently large.  Combining our results so far we conclude that $M_\gamma$ does not contain ITEs for $\gamma>0$, $n_\gamma$ chosen properly and all $\widetilde{n}^\ast\geq n_\gamma$ of the form \eqref{ratio}.
	
	\noindent Finally, we examine the open sectors $R_\gamma$ and $L_\gamma$ and prove that also these are ITE-free if $\widetilde{n}^\ast$ is sufficiently large and $\gamma>0$ chosen properly. The idea will be to show that any solution to \eqref{dkappa} located in $R_\gamma$ (in $L_\gamma$)  emerges from (arrives at) either the radial boundary part of $R_\gamma$ (of $L_\gamma$) or the real axis for some $n$, both of which will be impossible for ITE trajectories if $\gamma>0$ is small and $n$ large enough, cf. \Cref{figX}. Again, we only present the proof for $R_\gamma$ since $L_\gamma$ can be treated along the same lines. We start with validating the impossibility of intersection points of complex-valued ITE trajectories with $\partial R_\gamma$: We already know by \Cref{prop} that the real axis cannot be intersected apart from $\kappa^\ast$ for any $n\ne 1$ if $\gamma>0$ is small enough such that $\kappa^\ast$ is the only root of $J_p$ in $B_\gamma(\kappa^\ast)$. Further, ITEs $\kappa_n$ with $|\kappa_n|=\gamma$ and $0<\arg(\kappa_n-\kappa^\ast)<\pi/3$ cannot exist either for $n$ large, say $n\geq n_\partial$, which can be shown exactly as for the sectors \eqref{sectors} above. Hence, we are left to prove that \eqref{dkappa} forces any solution located in $R_\gamma$ for $n\geq n_\gamma>n_\partial$ large enough to origin from one of the two inadmissible fractions of $\partial R_\gamma$ excluding $\kappa^\ast$ for some $n\geq n_\partial$. For this purpose, we define as in the proof of \Cref{angle} for $\kappa \in R_\gamma$ and arbitrary $n\ne 1$
	\begin{align*}
		d_n(\kappa):=-\frac{n\kappa^2-p^2}{2n(n-1)\kappa}-\frac{g_p(\kappa)}{2n(n-1)(\kappa-\kappa^\ast)^2}
		=\frac{-\kappa+\frac{p^2}{n\kappa}-\frac{g_p(\kappa)}{n(\kappa-\kappa^\ast)^2}}{2(n-1)}\ ,
	\end{align*}
	where again
	\begin{align*}
		g_p(\kappa):=\frac{\kappa (\kappa-\kappa^\ast)^2 J'_p(\kappa)^2}{J_p(\kappa)^2}\ .
	\end{align*}
	Then, we use an orthogonal decomposition in the complex plane of the form
	\begin{align*}
		(n-1)d_n(\kappa)=a_nv+b_n\mathrm{e}^{\mathrm{-i\frac{\pi}{2}}}v
	\end{align*}
	with $a_n,b_n\in \mathbb{R}$ and will prove that there is $c>0$ such that 
	\begin{align}\label{sup}
		\sup_{(n_c,\infty)}b_n<-c 
	\end{align}
	independently of $\kappa\in R_\gamma$ if $n_c>1$ is chosen appropriately. 
	Noting that the distance from any $\kappa\in R_\gamma$ in the direction of $\mathrm{e}^{\mathrm{-i\frac{\pi}{2}}}v$ to the real axis or to the radial boundary part of $R_\gamma$ is bounded by $\cos(\pi/3)\gamma$, we can find $\Delta n>0$ such that 
	\begin{align}\label{distance}
		\int_{n_c}^{n_c+\Delta n}\frac{c}{n-1}\,\mathrm{d}n>\cos(\pi/3)\gamma\ .
	\end{align}
	If there was some ITE $\kappa_n \in R_\gamma$ for $n\geq n_c+\Delta n$, \eqref{distance} would enforce an intersection of the corresponding ITE trajectory with either the radial boundary part of $\partial R_\gamma$ or with the real axis apart of $\kappa^\ast$ for some refractive index greater than $n_c$, where the $\kappa^\ast$-exclusion follows by the fact that $c>0$ in \eqref{sup} for all $\kappa\in R_\gamma$. Choosing $n_c\geq n_\partial$, we recall that such intersections do not exist and may conclude that $R_\gamma$ is ITE-free for $n\geq n_\gamma\geq n_c+\Delta n$. 
	
	\noindent In order to prove \eqref{sup}, we note that if $g_p$ within the definition of $d_n$ was positive for all $\kappa\in R_\gamma$, existence of $c>0$ would follow immediately because the $(\mathrm{e}^{\mathrm{-i\frac{\pi}{2}}}v)$-component of both $-\kappa+p^2/(n\kappa)$ and $-g_p(\kappa)/(n(\kappa-\kappa^\ast)^2)$ are then strictly negative and non-positive, respectively, provided $n>1$ and $\gamma>0$ is sufficiently small. According to the proof of \Cref{angle} though, see especially \eqref{gp2nd} and its paragraph above, we actually have $\theta_\gamma<\arg(g_p(\kappa))\leq 0$ modulo $2\pi$ for $\kappa\in B_\gamma(\kappa^\ast)$ with $\theta_\gamma\to 0$ as $\gamma\to 0$. Therefore, $-g_p(\kappa)/(n(\kappa-\kappa^\ast)^2)$ rather has a positive  $(\mathrm{e}^{\mathrm{-i\frac{\pi}{2}}}v)$-component for $\pi/3-2\theta_\gamma\leq\arg(\kappa-\kappa^\ast)<\pi/3$ modulo $2\pi$ and $|\kappa-\kappa^\ast|<\gamma$ 
	which is given by 
	\begin{align}\label{drehen}
		\begin{split}
			&\left|\frac{g_p(\kappa)}{n(\kappa-\kappa^\ast)^2}\right|\tan\left(\frac{\pi}{3}-\arg\left(-\frac{g_p(\kappa)}{n(\kappa-\kappa^\ast)^2}\right)\right)\\
			=&\left|\frac{g_p(\kappa)}{n(\kappa-\kappa^\ast)^2}\right|\tan\left(\frac{\pi}{3}-\arg\big(-(\kappa-\kappa^\ast)^{-2}\big)-\arg\big(g_p(\kappa)\big)\right)
		\end{split}
	\end{align}
	and which we aim to compensate by adding $-\kappa+p^2/(n\kappa)$ for $n$ large.
	Since $g_p'(\kappa^\ast)=0$ and $g_p''(\kappa^\ast)<0$, cf. \eqref{gp2nd} and its paragraph above, the second order approximation of $g_p$ around $\kappa^\ast$ ensures that for $\gamma>0$ small there exists $C_\gamma>0$ such that
	\begin{align*}
		0\geq\arg(g_p(\kappa))\geq -\tan\left(\frac{C_\gamma|\kappa-\kappa^\ast|^2}{\kappa^\ast}\right)\quad \text{mod }2\pi
	\end{align*}
	for all $\kappa\in R_\gamma$. Since $\arg\big(-(\kappa-\kappa^\ast)^{-2}\big)\geq \pi/3$ modulo $2\pi$ for $\kappa\in R_\gamma$, we can estimate
	\begin{align}\label{term1}
		\left|\frac{g_p(\kappa)}{n(\kappa-\kappa^\ast)^2}\right|\tan\left(\frac{\pi}{3}-\arg\left(-\frac{g_p(\kappa)}{n(\kappa-\kappa^\ast)^2}\right)\right)\leq \left|\frac{g_p(\kappa)}{n(\kappa-\kappa^\ast)^2}\right|\tan\left(\frac{C_\gamma|\kappa-\kappa^\ast|^2}{\kappa^\ast}\right)\leq \frac{4C_\gamma}{n}\ ,
	\end{align}
	where $\gamma>0$ is additionally restricted to fulfill $\tan(t)\leq 2t$ for all $0\leq t<\gamma$ and $|g_p(\kappa)|\leq 2\kappa^\ast$ for all $\kappa\in R_\gamma$. In contrast, the $(\mathrm{e}^{\mathrm{-i\frac{\pi}{2}}}v)$-component of $-\kappa+p^2/(n\kappa)$ reads 
	\begin{align}\label{term2}
		\left|-\kappa+\frac{p^2}{n\kappa}\right|\tan\left(\arg\left(-\kappa+\frac{p^2}{n\kappa}\right)-\frac{\pi}{3}\right)\leq  \tan\left(\frac{2\pi}{3}\right)\frac{\kappa^\ast}{2}<0
	\end{align}
	which is independent of $n$ for all $\kappa\in R_\gamma$ with $\gamma>0$ small and $n$ sufficiently large. In particular, adding the left-hand side of \eqref{term1} and \eqref{term2} each, we deduce that there is indeed $c>0$ such that \eqref{sup} holds. Hence, $R_\gamma$ (and likewise $L_\gamma$) is ITE-free for $\gamma>0$, $n_\gamma$ chosen properly and all $n\geq n_\gamma$. The proof is now complete.
\end{proof}

\begin{proposition}\label{prop}
	Let $p\in\mathbb{N}_0$ and $\kappa^\ast$ be such that $J_p(\kappa^\ast)=0$. Further, assume that $\phi\in\{\pi/3,2\pi/3\}$ is fixed in $v=\mathrm{e}^{\mathrm{i}\phi}$ and let $g_1$ be the unique solution to
	\begin{align}\label{ode}
		g_n'(t)=v\left(1+\frac{1}{\kappa^\ast+tv}g_n(t)+n\frac{(\kappa^\ast+tv)^2-\frac{p^2}{n}}{(\kappa^\ast+tv)^2}g_n(t)^2\right)\ ,\quad g_n(0)=0\ .
	\end{align}
	Then we can find $\gamma>0$ and $n_\gamma>1$ such that for any solution $g_n:[0,\gamma]\to\mathbb{C}$ to \eqref{ode} with $n\geq n_\gamma$ we have that
	\begin{align*}
		|g_n(t)|<|g_1(t)|\qquad \text{for all }\ 0<t\leq \gamma\ .
	\end{align*}
\end{proposition}
\begin{proof}
	Uniqueness of solutions to \eqref{ode} follows by the Picard-Lindelöf theorem and existence of even global solutions $g_n:[0,\infty)\to\mathbb{C}$, that are solutions without blow-up for some $t<\infty$, is ensured by \eqref{gn} for infinitely many $n\nearrow \infty$. Since \eqref{ode} implies \eqref{taylor1} and \eqref{taylor2} even for all $n\geq 1$, the proof of the previous lemma can be adopted and yields that for each $n$ there is $\varepsilon>0$ (depending on $n$) such that $|g_n(t)|<|g_1(t)|$ for all $0<t<\varepsilon$. In order to replace $\varepsilon>0$ by some $n$-independent threshold $\gamma>0$, we assume contrarily that for any $\gamma>0$ sufficiently small we can find a sequence of $n$ tending to infinity such that
	\begin{align}\label{tn}
		t_n:=\inf\{0<t\leq \gamma:\ |g_n(t)|\geq|g_1(t)|\}
	\end{align} 
	is well-defined each. In particular, it holds that $t_n\geq \varepsilon>0$ and by continuity
	\begin{align*}
		|g_n(t_n)|=|g_1(t_n)|\ .
	\end{align*}
	Next, we compute with the help of \eqref{ode} for $n\geq 1$
	\begin{align}\label{ableit}
		\left[\frac{|g_n(t)|^2}{2}\right]'=\Re\left(\overline{v}g_n(t)\right)+\Re\left(\frac{v}{\kappa^\ast+tv}\right)|g_n(t)|^2+n\Re\left(v\frac{g_n(t)}{|g_n(t)|}\frac{(\kappa^\ast+tv)^2-\frac{p^2}{n}}{(\kappa^\ast+tv)^2}\right)|g_n(t)|^3\ .
	\end{align}
	Since $\kappa^\ast>p$ according to \cite[p. 486]{Wa66}, we restrict to $\gamma>0$ small such that 
	\begin{align}\label{constraint}
		0\leq \arg\left(\frac{(\kappa^\ast+tv)^2-\frac{p^2}{n}}{(\kappa^\ast+tv)^2}\right)\leq\arg\left(\frac{(\kappa^\ast+tv)^2-p^2}{(\kappa^\ast+tv)^2}\right)<\frac{\pi}{12}\quad \text{mod }2\pi
	\end{align}
	for all $0<t\leq \gamma$ and $n\geq 1$. Further, we assume that $\gamma$ is such that $|g_1|$ is monotonically increasing in $(0,\gamma)$ with $t\geq |g_1(t)|>t/2$ and
	\begin{align}\label{low}
		\begin{split}
			\frac{\pi}{3}\leq \arg(g_1(t)), \arg(g_1'(t)) \leq\frac{5\pi}{12}\quad &\text{mod }2\pi,\quad \text{if }v=\mathrm{e}^{\mathrm{i}\frac{\pi}{3}}\ ,\\
			\frac{7\pi}{12}\leq \arg(g_1(t)), \arg(g_1'(t))\leq \frac{2\pi}{3}\quad &\text{mod }2\pi,\quad \text{if }v=\mathrm{e}^{\mathrm{i}\frac{2\pi}{3}}\ ,
		\end{split}
	\end{align}
	which is feasible due to \eqref{taylor1}, and we impose that
	\begin{align*}
		\left|\arg\left(1+\frac{1}{\kappa^\ast+tv}g_n(t)\right)\right|<
		\frac{\pi}{12}\quad \text{mod }2\pi
	\end{align*}
	for all $0<t<t_n$ which is justified by exploiting $|g_n(t)|<|g_1(t)|\leq |g_1(t_n)|\leq |g_1(\gamma)|$, cf. \eqref{tn}. In what follows we thus have to distinguish the cases $\phi=\pi/3$ and $\phi=2\pi/3$, which can be treated similarly though, so we only present our proof for $v=\mathrm{e}^{\mathrm{i}\frac{\pi}{3}}$ in the sequel. We can then conclude that $g_n(t)\ne 0$ and
	\begin{align}\label{sector}
		\frac{\pi}{4}\leq \arg(g_n(t))\leq \frac{3\pi}{4}\quad \text{mod }2\pi
	\end{align}
	for all $0<t<t_n$ since for $\arg(g_n(t))=\pi/4$ modulo $2\pi$ or $\arg(g_n(t))=3\pi/4$ modulo $2\pi$ we have by checking \eqref{ode} summandwise that $\pi/4<\arg(g_n'(t))< 11\pi/12$ modulo $2\pi$ or $-\pi/4<\arg(g_n'(t))< 5\pi/12$ modulo $2\pi$, respectively, that is the normal component of $g_n'(t)$ relative to $g_n(t)$ points into the sector \eqref{sector}. For $g_n(t)=0$ we similarly get $g_n'(t)=v=\mathrm{e}^{\mathrm{i}\frac{\pi}{3}}$, so that if $g_n(t)$ once lies in the sector it remains inside for all $t<t_n$. 
	To proceed, we will have to consider the cases $\arg(g_n(t_n))\geq \arg(g_1(t_n))$ modulo $2\pi$ and $\arg(g_n(t_n))< \arg(g_1(t_n))$ modulo $2\pi$ separately:
	
	\noindent If $\arg(g_n(t_n))\geq \arg(g_1(t_n))$ modulo $2\pi$, \eqref{low} and \eqref{sector} together yield that
	\begin{align}\label{angbounds}
		\frac{\pi}{3}\leq \arg(g_1(t_n)) \leq \arg(g_n(t_n))\leq \frac{3\pi}{4}\quad \text{mod }2\pi\ ,
	\end{align}
	so \eqref{ableit} can be estimated by
	\begin{align*}
		\left[\frac{|g_n|^2}{2}\right]'_{|t=t_n}=&\ \Re\left(\overline{v}\frac{g_n(t_n)}{|g_n(t_n)|}\right)|g_n(t_n)|+\Re\left(\frac{v}{\kappa^\ast+t_nv}\right)|g_n(t_n)|^2\\
		&+n\Re\left(v\frac{g_n(t)}{|g_n(t)|}\frac{(\kappa^\ast+t_nv)^2-\frac{p^2}{n}}{(\kappa^\ast+t_nv)^2}\right)|g_n(t_n)|^3\\
		\leq&\  \Re\left(\overline{v}\frac{g_1(t_n)}{|g_1(t_n)|}\right)|g_1(t_n)|+\Re\left(\frac{v}{\kappa^\ast+t_nv}\right)|g_1(t_n)|^2\\
		&+n\Re\left(v\frac{g_n(t)}{|g_n(t)|}\frac{(\kappa^\ast+t_nv)^2-\frac{p^2}{n}}{(\kappa^\ast+t_nv)^2}\right)|g_1(t_n)|^3\ .
	\end{align*}
	Using \eqref{constraint} and the angular bounds in \eqref{angbounds}, we can find $n_\gamma> 1$ such that 
	\begin{align*}
		n\Re\left(v\frac{g_n(t_n)}{|g_n(t_n)|}\frac{(\kappa^\ast+t_nv)^2-\frac{p^2}{n}}{(\kappa^\ast+t_nv)^2}\right)<\Re\left(v\frac{g_1(t_n)}{|g_1(t_n)|}\frac{(\kappa^\ast+t_nv)^2-p^2}{(\kappa^\ast+t_nv)^2}\right)<0
	\end{align*}
	for all $n\geq n_\gamma$. Hence, we obtain $\left[|g_n|^2/2\right]'_{|t=t_n}<\left[|g_1|^2/2\right]'_{|t=t_n}$, which contradicts the definition of $t_n$ in \eqref{tn}.
	
	\noindent In the other case, that is $\arg(g_n(t_n))< \arg(g_1(t_n))$ modulo $2\pi$, there exist $0<\widetilde{s}_n<\widetilde{t}_n<t_n$ such that $g_1(\widetilde{s}_n)=g_n(\widetilde{t}_n)$ and 
	\begin{align}\label{left}
		\arg(-g_1'(\widetilde{s}_n))\leq
		\arg(g_n'(\widetilde{t}_n))\leq \arg(g_1'(\widetilde{s}_n))\quad \text{mod }2\pi\ .
	\end{align}
	More precisely, the existence of $\widetilde{s}_n,\widetilde{t}_n$ follows by our particular assumption and $\arg(g_n(t))\geq \arg(g_1(t))$ in a neighborhood of $t=0$ which is a direct consequence of $g_1'(0)=g_n'(0)=v=\mathrm{e}^{\mathrm{i}\frac{\pi}{3}}$, $g_n'''(0)<g_1'''(0)<0$ and
	\begin{align*}
		g_{n}(t)-g_1(t)=\int_0^t \big(g_{n}'''(s)-g_1'''(s)\big)\frac{(t-s)^2}{2}\,\mathrm{d}s\ .
	\end{align*}
	The constraint \eqref{left} can be restated as
	\begin{align}\label{constraint2}
		\arg(-g_1'(\widetilde{s}_n))\leq\arg(g_n'(\widetilde{t}_n)-g_1'(\widetilde{s}_n))\leq \arg(g_1'(\widetilde{s}_n))\quad \text{mod }2\pi\ .
	\end{align}
	Using \eqref{ode} and $g_1(\widetilde{s}_n)=g_n(\widetilde{t}_n)$, we can write
	\begin{align}
		\begin{split}\label{upper}
			g_n'(\widetilde{t}_n)-g_1'(\widetilde{s}_n)=&\ 
			\frac{\widetilde{s}_n-\widetilde{t}_n}{(\kappa^\ast+\widetilde{s}_nv)(\kappa^\ast+\widetilde{t}_nv)}v^2g_1(\widetilde{s}_n)+
			(n-1)vg_1(\widetilde{s}_n)^2\\
			&- \frac{p^2(\widetilde{s}_n-\widetilde{t}_n)\big(2\kappa^\ast+(\widetilde{t}_n+\widetilde{s}_n)v\big)}{(\kappa^\ast+\widetilde{s}_nv)^2(\kappa^\ast+\widetilde{t}_nv)^2}v^2g_1(\widetilde{s}_n)^2\\
			=&\ 
			(\widetilde{s}_n-\widetilde{t}_n)g_1(\widetilde{s}_n)v^2 \frac{(\kappa^\ast+\widetilde{s}_nv)(\kappa^\ast+\widetilde{t}_nv)-p^2g_1(\widetilde{s}_n)\big(2\kappa^\ast+(\widetilde{t}_n+\widetilde{s}_n)v\big)}{(\kappa^\ast+\widetilde{s}_nv)^2(\kappa^\ast+\widetilde{t}_nv)^2}\\
			&+(n-1)vg_1(\widetilde{s}_n)^2
			\ .
		\end{split}
	\end{align}
	Additionally, we can assume that $\gamma>0$ is small enough such that 
	\begin{align*}
		\left|\arg\left(\frac{(\kappa^\ast+\widetilde{s}_nv)(\kappa^\ast+\widetilde{t}_nv)-p^2g_1(\widetilde{s}_n)\big(2\kappa^\ast+(\widetilde{t}_n+\widetilde{s}_n)v\big)}{(\kappa^\ast+\widetilde{s}_nv)^2(\kappa^\ast+\widetilde{t}_nv)^2}\right)\right|<\frac{\pi}{12}\quad \text{mod }2\pi
	\end{align*}
	and
	\begin{align*}
		\left|\frac{(\kappa^\ast+\widetilde{s}_nv)(\kappa^\ast+\widetilde{t}_nv)-p^2g_1(\widetilde{s}_n)\big(2\kappa^\ast+(\widetilde{t}_n+\widetilde{s}_n)v\big)}{(\kappa^\ast+\widetilde{s}_nv)^2(\kappa^\ast+\widetilde{t}_nv)^2}\right|<\frac{2}{{\kappa^\ast}^2}
	\end{align*}
	for all $0<\widetilde{s}_n<\widetilde{t}_n<\gamma$. Because of $\widetilde{s}_n\geq|g_1(\widetilde{s}_n)|>\widetilde{s}_n/2$ and $\pi/3\leq\arg(g_1(\widetilde{s}_n)),  \arg(g_1'(\widetilde{s}_n)) \leq 5\pi/12$ by \eqref{low}, we conclude with \eqref{upper}, in order for \eqref{constraint2} to hold as $n\to\infty$, that there is $c_\gamma>0$ such that 
	\begin{align*}
		\frac{\widetilde{t}_n}{n\widetilde{s}_n}\geq c_\gamma
	\end{align*}
	for all $n\geq n_{\gamma}$, where $n_{\gamma}>1$ is  sufficiently large. Hence, we get that
	\begin{align}\label{contra1}
		|g_n(\widetilde{t}_n)|=|g_1(\widetilde{s}_n)|\leq \widetilde{s}_n\leq \frac{\widetilde{t}_n}{c_\gamma n}
	\end{align}
	for all $n\geq n_{\gamma}$. We can even assume that $g_n$ is monotonically increasing in $(0,\widetilde{t}_n)$ for $n\geq n_{\gamma}$ since $\left[|g_n|^2/2\right]'\leq 0$ requires by \eqref{ableit} that
	\begin{align}\label{lowerbnd}
		|g_n(t)|\geq \sqrt{\frac{\Re\left(\overline{v}\frac{g_n(t)}{|g_n(t)|}\right)}{n\Re\left(-v\frac{g_n(t)}{|g_n(t)|}\frac{(\kappa^\ast+t_nv)^2-\frac{p^2}{n}}{(\kappa^\ast+t_nv)^2}\right)}}\geq \frac{d_\gamma}{\sqrt{n}}
	\end{align}
	for some $d_\gamma>0$ which is independent of $n$ for $0<t<t_n$ due to \eqref{sector} and \eqref{constraint}. So if $|g_n(\widetilde{h}_n)|>|g_n(\widetilde{t}_n)|$ for some $\widetilde{h}_n<\widetilde{t}_n$, and since also $|g_n(t_n)|=|g_1(t_n)|\geq |g_1(\widetilde{s}_n)|\geq |g_n(\widetilde{t}_n)|$, $|g_n|$ would have a local minimum $\widetilde{i}_n\in(\widetilde{h}_n,t_n)$ with $|g_n(\widetilde{i}_n)|\leq |g_n(\widetilde{t}_n)|\leq \gamma/(c_\gamma n)$ by \eqref{contra1}, which contradicts \eqref{lowerbnd} for $t=\widetilde{i}_n$ and $n$ large.
	Hence, we have that $\max_{[0,\widetilde{t}_n]}|g_n|\leq \gamma/(c_\gamma n)$ for $n\geq n_{\gamma}$ large and therefore
	\begin{align*}
		|g_n(\widetilde{t}_n)|&=\left|\int_0^{\widetilde{t}_n}g_n'(t)\,\mathrm{d}t\right| \\
		&=\left|\int_0^{\widetilde{t}_n}v\left(1+\frac{1}{\kappa^\ast+tv}g_n(t)+n\frac{(\kappa^\ast+tv)^2-\frac{p^2}{n}}{(\kappa^\ast+tv)^2}g_n(t)^2\right)\,\mathrm{d}t\right|\\
		&\geq \widetilde{t}_n-\int_0^{\widetilde{t}_n}\frac{|g_n(t)|}{\kappa^\ast}+n|g_n(t)|^2\,\mathrm{d}t\\
		&\geq \widetilde{t}_n\left(1-\frac{\gamma}{\kappa^\ast c_\gamma n}-\frac{\gamma^2}{c_\gamma^2 n}\right)\ .
	\end{align*}
	This is a contradiction to \eqref{contra1} for $n\geq n_{\gamma}$ large and thus to \eqref{tn} again, which completes the proof.
\end{proof}

\begin{remark}\label{3D}
	To extend our findings from the unit disk to the unit ball in 3D, the Fourier Bessel ansatz \Cref{eigf} in 2D for ITP eigenfunctions needs to be replaced by
	\begin{align*}
		\begin{split}
			v_n(r,\phi,\theta)&=j_p(\kappa_nr)\cos(p\phi)P_p^\ell(\cos(\theta))\ , \\ \Big(v_n(r,\phi,\theta)&=j_p(\kappa_nr)\sin(p\phi)P_p^\ell(\cos(\theta))\ , p\ne 0\Big)
		\end{split}
	\end{align*}
	and
	\begin{align*}
		\begin{split}
			w_n(r,\phi,\theta)&=\alpha_n j_p(\sqrt{n}\kappa_nr)\cos(p\phi)P_p^\ell(\cos(\theta))\ ,\\ \Big(w_n(r,\phi,\theta)&=\alpha_nj_p(\sqrt{n}\kappa_nr)\sin(p\phi)P_p^\ell(\cos(\theta))\ , p\ne 0\Big),
		\end{split}
	\end{align*}
	where $p\in \mathbb{N}_0$ and $\ell\in\{-p,\ldots,p\}$. Here, $P_p^\ell$ is the associated Legendre polynomial and $j_p(z)=\sqrt{\frac{\pi}{2z}}J_{p+\frac{1}{2}}(z)$ is the spherical Bessel function of the first kind of order $p$ satisfying the second-order ordinary differential equation
	\begin{align}\label{sbes}
		z^2j''_p(z)+2z j'_p(z)+\left(z^2-p(p+1)\right)=0\ .
	\end{align} 
	The characteristic equation \Cref{determinant} still keeps the same structure, that is
	\begin{align}\label{determinant2}
		f_p(n,\kappa):=\kappa j_p'(\kappa)j_p(\kappa \sqrt{n})-\kappa \sqrt{n}j_p(\kappa)j_p'(\kappa \sqrt{n})\ .
	\end{align}
	We will present numerical studies below which demonstrate the similarities of corresponding ITE trajectories.
\end{remark}

\section{Numerical observations}\label{sec3}
In this section, we present numerical results for some standard scattering shapes to visualize our theoretical findings and to illustrate further interesting phenomena. The underlying Matlab program can be downloaded from:\\ \texttt{https://github.com/kleefeld80/ITEtrajectory} .

\begin{figure}[htbp]
	\centering
	\includegraphics[width=10cm]{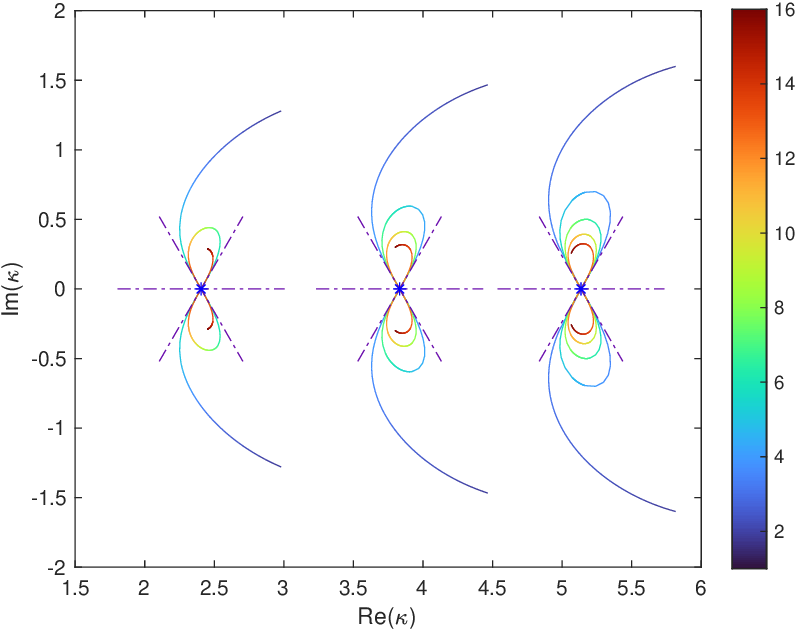}
	\caption{\label{fig1}The first three complex-conjugated pairs of  complex-valued ITE trajectories for the unit disk and $p\in\{0,1,2\}$ without counting multiplicity, respectively,  using $n\in (1,16]$.}
\end{figure}

\subsection{The unit disk for $n>1$}
The first complex-valued roots $\kappa_n$ of the function $F_p$ from \Cref{determinant} (ordered by real parts) are computed for $p\in\{0,1,2\}$ and sufficiently many $n\in (1,16]$ via Beyn's second integral algorithm, see \cite[p. 3860]{beyn2012integral}. Using interpolation with respect to $n$ to obtain smooth trajectories, \Cref{fig1} shows the resulting output where the color bar refers to the varying index of refraction $n>1$.

As we observe and expect, the ITE trajectories come in complex-conjugated pairs and arise at certain points satisfying $\mathrm{Im}(\kappa_n)\neq 0$ for $n$ close to one. For increasing $n$ they reapproach the first three DELs of the unit disk, cf. \Cref{theorem210}, which are marked by blue asterisks in \Cref{fig1} and are approximately given by $2.4048$, $3.8317$, and $5.1356$. For instance, the first complex-conjugated trajectory pair corresponds to $p=0$ and  passes through the approximate DEL $2.4048$ for the first time at $n\approx 5.2689$ and then for $n\approx 12.9491$. These $n$ coincide, in agreement with \Cref{ratio}, with the squared ratio of successive larger roots of $J_0$ than $2.4048$ and $2.4048$ itself, respectively. Recall that incident and outgoing angles at DELs are restricted to $\pm \pi/3$ by \Cref{angle} which is highlighted additionally in the figure by green lines. Further, we see that the trajectory pair converges as a whole towards $2.4048$ for growing $n$, cf. \Cref{theorem210}. The other two complex-conjugated trajectory pairs in \Cref{fig1} corresponding to $p=1$ and $p=2$ behave likewise and are not discussed in further detail.
Also note that we cannot find any other ITE trajectories than the DEL-recurrent ones in the complex plane which also follows by the uniqueness statement of \Cref{theorem210}.
Finally, we point out that although ITEs are formally not defined for $n=1$, \Cref{fig1} shows that complex-valued ITE trajectories converge each to a unique non-real number as $n\to 1$. We can even compute $\kappa_1:=\lim_{n\to 1}\kappa_n$ according to the following reasoning: for the latter limit to exist, $\kappa_n'$ must be locally integrable around $n=1$ by the fundamental theorem of calculus. Rewriting \eqref{dkappa} as
\begin{align*}
	\kappa_n'=-\frac{J_p(\kappa_n)^2(n\kappa_n^2-p^2)-\kappa_n^2 J'_p(\kappa_n)^2}{2n(n-1)\kappa_nJ_p(\kappa_n)^2}\ ,
\end{align*}
we conclude, due to the fact that $1/(n-1)$ is not integrable around $n=1$ locally, that the enumerator must vanish for $n\to 1$. Exploiting the modified Beyn algorithm again for finding the enumerator's roots, $\kappa_1$ for the 3 complex-conjugated pairs of complex-valued trajectories in \Cref{fig1} are approximately  $2.9804 \pm 1.2796\mathrm{i}$, $4.4663 \pm 1.4675\mathrm{i}$, and $5.8169 \pm 1.6000\mathrm{i}$ for $p=0$, $p=1$, and $p=2$, respectively. 

\begin{figure}[htbp]
	\centering
	\includegraphics[width=10cm]{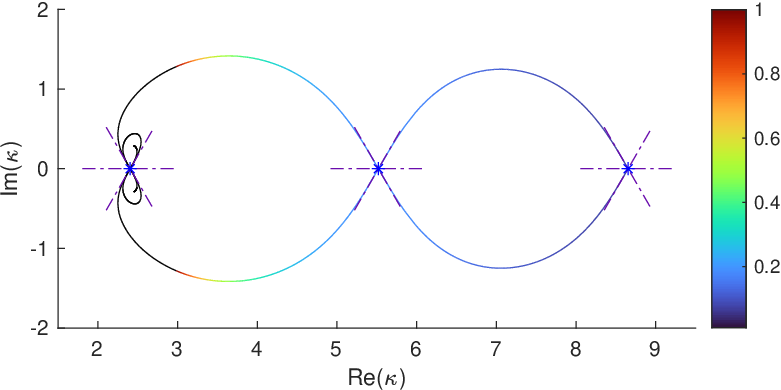}
	\caption{\label{fig2}The first complex-conjugated pair of  complex-valued ITE trajectories for the unit disk and $p=0$ from \Cref{fig1} 
		extended to $n\in (0,16]$ (color bar refers to $n<1$ only).}
\end{figure}

\subsection{The unit disk for $0<n<1$}
In \Cref{fig2}, we plot the trajectory continuation in $0<n<1$ for the first complex-conjugated pair of complex-valued ITE trajectories of the unit disk with Bessel index $p=0$ from the previous subsection. We observe that it connects continuously through the formal singularity $n=1$ to $n>1$ (shown here in black for $n\in (1,16]$ and copied from \Cref{fig1}). As before, we used Beyn's second integral algorithm to generate the trajectories.

We directly notice that the eigenvalue tajectories for $0<n<1$ behave different from the ones for $n>1$. Indeed, for $0<n<1$ the trajectories are not recurrent to a single DEL any more but escape to infinity as $n\to 0$, crossing the real axis at successive DELs instead. Also note that ingoing and outgoing directions at DELs are now given by $\pm 2\pi/3$, cf. \Cref{angle}. Still, the two regimes $0<n<1$ and $n>1$ are connected, even for arbitrary scatterers, via the relation $\kappa_n=\kappa_{1/n}/\sqrt{n}$ which follows by exchanging the roles of $v$ and $w$ in \eqref{ITP}. Thus, we restrict our numerical studies to $n>1$ in the sequel.

\begin{figure}[htbp]
	\centering
	\includegraphics[width=10cm]{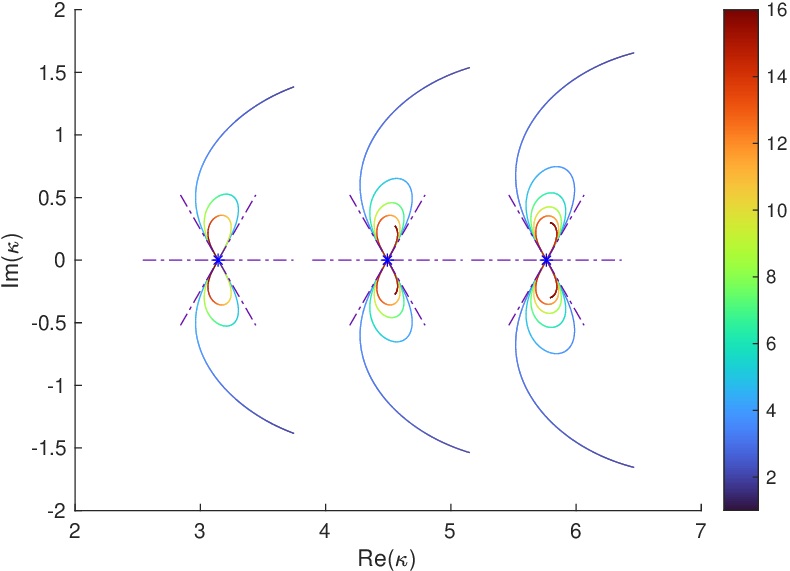}
	\caption{\label{fig3}The first three complex-conjugated pairs of  complex-valued ITE trajectories for the unit ball and $p\in\{0,1,2\}$ without counting multiplicity, respectively, using $n\in (1,16]$.}
\end{figure}

\subsection{The unit ball}
In \Cref{fig3}, we provide the first three complex-conjugated pairs of complex-valued ITE curves for the unit ball in 3D using $n\in (1,16]$ and $p\in \{0,1,2\}$, cf. \Cref{3D}. For this purpose, we need to compute the complex-valued roots of \Cref{determinant2} instead of \Cref{determinant} and can then proceed as in \Cref{fig1}.

We recognize an analogue behavior as for the unit disk. Specifically for $p=0$, which is the case of spherically symmetric ITP eigenfunctions, Colton \& Leung already pointed out in \cite{colton2013complex} that ITEs associated with spherical Bessel index $p=0$ can only be real-valued for $n=q^2$ or $n=1/q^2$, where $q\in\mathbb{N}$. From our perspective of ITE trajectories, this implies that there are  simultaneous intersections with infinitely many DELs for the same $n$. This can  be seen from $j_0(\kappa)=\sin(\kappa)/\kappa$ which has equidistant roots (and thus DELs) at $\kappa^\ast\in \mathbb{N}\pi$. Comparing with \Cref{ratio} yields $n^\ast=(qm\pi)^2/(m\pi)^2=q^2$ for any $m\in \mathbb{N}$, that is, complex-valued ITE trajectories intersect $\kappa^\ast=m\pi$ for all $m\in \mathbb{N}$ simultaneously whenever $n=q^2$, $q\in\mathbb{N}$. 
Such a simultaneous recurrence has not been observed for the unit disk since the roots of $J_0$, as well as of higher order Bessel functions, are only asymptotically equidistant, see \cite{abramowitz1964handbook,Be95}. 

\noindent Analogue to the unit disk, we can also define $\kappa_1$ as the limit of complex-valued ITE trajectories as $n\to 1$. For their computation in the 3D case, we use
\begin{align*}
	\kappa_n'= \frac{\big(n\kappa_n^2-p(p+1)\big)j_p(\kappa_n)^2+\kappa_nj_p(\kappa_n)j_p'(\kappa_n)+\kappa_n^2j_p'(\kappa_n)^2}{2n(n-1)\kappa_nj_p(\kappa_n)^2}\ ,
\end{align*}
cf. \eqref{sbes} and \eqref{determinant2}. Setting the latter enumerator to zero, the modified Beyn algorithm yields the solutions $3.7488 \pm 1.3843\mathrm{i}$, $5.1524 \pm  1.5380\mathrm{i}$, and $6.4652 \pm 1.6559\mathrm{i}$ using $p=0$, $p=1$, and $p=2$, respectively.

\begin{figure}[htbp]
	\centering
	\includegraphics[width=10cm]{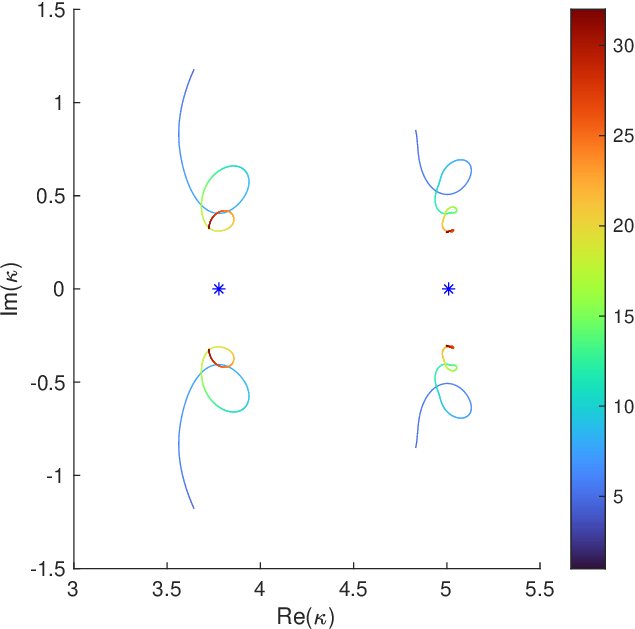}
	\caption{\label{fig4} The first two complex-conjugated pairs of  complex-valued ITE trajectories for the ellipse with semi-axes $1$ and $0.5$ using $n\in [4,32]$.}
\end{figure}

\subsection{The ellipse}
Next, we consider an ellipse with semi-axes $1$ and $0.5$ to also investigate non-spherical scatterers. Since a separation of variable ansatz does not simplify to a decoupled analytical expression as in \Cref{determinant} or \Cref{determinant2}, we employ the modified method of fundamental solution (modified MFS) to compute approximate ITEs, see \cite[Section 3.2]{kleefeldpieronek2018} for more details on this algorithm. In its original notation, we placed $m_I=10$ nodes on a circle with radius $0.4$ inside the ellipse, $m=40$ collocation points along the ellipse's boundary and $m=40$ source points on an exterior circle with radius $4$. When refering to the modified MFS in the remainder of this section, we assume that all the auxiliary circles have the same center as the scattering object itself. Further, since ITEs are computed within the modified MFS as minimizers of some boundary collocation misfit function, we take the computed ITE from the previous $n$ iteratively as initial guess for the next when incrementing $n$. Thus we are only left to set one independent initial guess for the complex-valued ITE at the minimal $n$ of interest which we will fix as $n=4$. 

For our ellipse, we picked $4\pm \mathrm{i}$ and $5\pm \mathrm{i}$ as independent initial guesses to compute the first two complex-conjugated pairs of complex-valued ITE trajectories, respectively. The resulting curves for $n\in [4,32]$ are shown in \Cref{fig4}. Unlike for the disk or the ball, however, the ITE trajectories are not recurrent with respect to DELs any more, which are computed approximately as $3.7771$ and $5.0102$ according to the formula given in \cite[pp. 9]{grebenkov}. However, they tend to spiral down towards a unique DEL, cf. \Cref{big2}, without touching the real axis at all. Also, we still observe a one-to-one correspondence between DELs and complex-valued ITE trajectories governed by their seemingly convergent behavior as $n\to \infty$.

\begin{figure}[htbp]
	\centering
	\includegraphics[width=10cm]{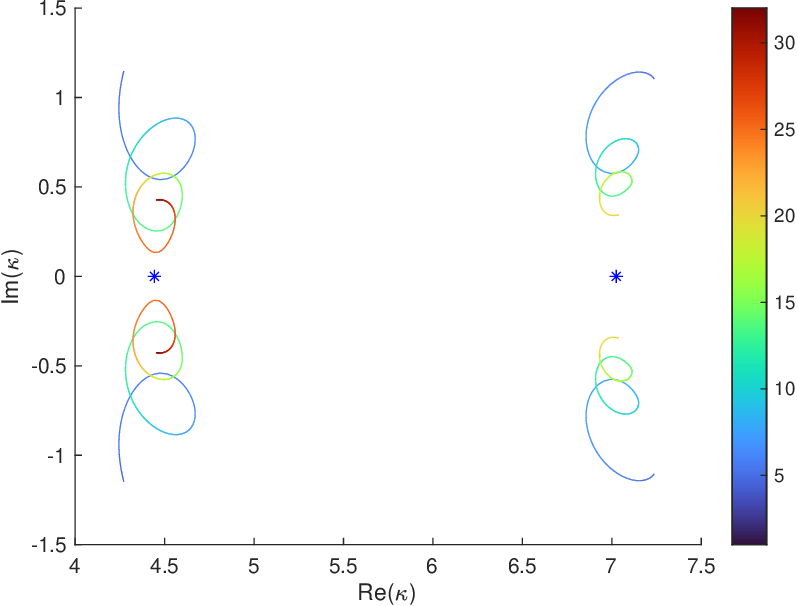}
	\caption{\label{fig5}The first two complex-conjugated pairs of  complex-valued ITE trajectories for the unit square using $n\in [4,32]$ and $n\in [4,20]$, respectively.}
\end{figure}

\subsection{The unit square}
In order to include non-smooth scatterers into our scope of investigation, we turn our attention now to the unit square. Here, the modified MFS has been exploited with $m_I=20$ nodes on an interior circle with radius $0.25$, $m=61$ collocation points along the boundary of the square apart from its corners (as the outer normal $\nu$ would not be defined otherwise) and $m=61$ source points on an exterior circle with radius $0.75$. The first two complex-conjugated pairs of complex-valued ITE trajectories are displayed in \Cref{fig5}, taking $4.5\pm\mathrm{i}$ and $7\pm\mathrm{i}$ as independent initial guesses at $n=4$, respectively. 

In contrast to the ellipse, the two trajectory pairs are now plotted with respect to different yet overlapping parameter domains of $n$ which is $[4,32]$ for the first and $[4,20]$ for the second. The reason is that the modified MFS generally suffers from ill-conditioning effects for large wave numbers. In particular, neither $\kappa_n$ in $v_n$ nor $\sqrt{n}\kappa_n$ in $w_n$ can be too large within \Cref{eigf}, so that especially the latter restricts the feasible choices of $n$.
In our plot both resulting trajectory pairs show again a spiral pattern which approach for growing $n$ the separated blue asterisks on the real axis, respectively, as would be expected by \Cref{big2}. They are given by $\sqrt{2}\pi\approx 4.4429$ and $\sqrt{5}\pi\approx 7.0248$ according to the first two DELs of the unit square in \cite[pp. 6]{grebenkov}. We still observe a one-to-one correspondence between DELs and complex-conjugated pairs of complex-valued ITE trajectories similar to the ellipse.


\begin{figure}[htbp]
	\centering
	\includegraphics[width=10cm]{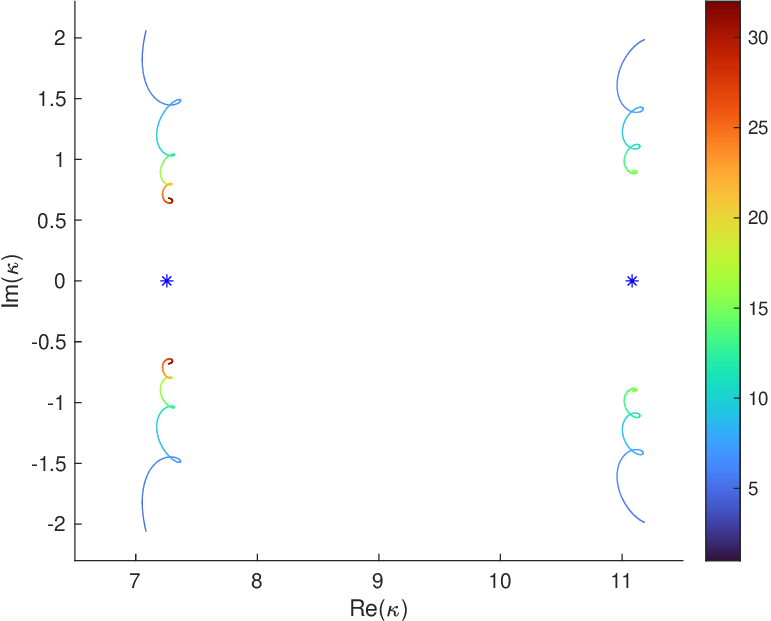}
	\caption{\label{fig6}The first two complex-conjugated pairs of  complex-valued ITE trajectories for the equilateral triangle with side length one using $n\in [4,32]$ and $n\in [4,16]$, respectively.}
\end{figure}

\subsection{The equilateral triangle}
We also consider an equilateral triangle with side length one for which the first two DELs are given by $4\pi/\sqrt{3}\approx 7.255$ and $4\pi\sqrt{7}/3\approx 11.082$, see \cite[pp. 10--11]{grebenkov}). Using the same parameters within the modified MFS as for the unit square 
but reducing $m=61$ to $m=51$ for the first trajecory pair 
yields the two complex-conjugated pairs of complex-valued ITE trajectories shown in \Cref{fig6}. The domain of $n$ was $n\in [4,32]$ for the first pair and $n\in [4,16]$ for the second. As independent initial ITE guesses at $n=4$ we took $7.3\pm 1.5\mathrm{i}$ and $11\pm 2\mathrm{i}$, respectively.
Altogether, we again observe a one-to-one correspondence between DELs and complex-conjugated pairs of complex-valued ITE trajectories. 


\subsection{The deformed ellipse}
We also present an example of a non-convex scatterer which is parameterized for $t\in [0,2\pi)$ by 
\begin{align}\label{deform}
	t \mapsto 
	\begin{pmatrix}
		0.75\cos(t)+0.3 \cos(2t)\\
		\sin(t)
	\end{pmatrix}\ .
\end{align}
Its exact shape is illustrated in \Cref{tikz1}.
\begin{figure}[htbp]
	\begin{center}
		\begin{tikzpicture}
			\begin{axis}[scale=0.5,
				trig format plots=rad,
				axis equal,
				hide axis,
				legend columns=4, 
				legend style={
					legend image post style={thick},
					nodes={scale=0.75},
					at={(0.5,0.2)},anchor=south,
				},
				]
				\addplot [domain=0:2*pi, samples=200, gray,opacity=1.0,fill=gray,fill opacity=0.4] ({6.0+0.75*cos(x)+0.3*cos(2*x)}, {sin(x) });
				\addplot [domain=0:2*pi, samples=200, black] ({6.0+0.75*cos(x)+0.3*cos(2*x)}, {sin(x) });
			\end{axis}
		\end{tikzpicture}
	\end{center}
	\caption{The deformed ellipse parameterized by \Cref{deform}
	}\label{tikz1}
\end{figure}
In this case, we employ the modified MFS both for computing the first two DELs of the deformed ellipse and for generating the first two pairs of complex conjugated ITE trajectories. Specifically, we distributed $m_I=20$ nodes on an interior circle with radius $0.2$ centered in the origin, $m=51$ collocation points along the boundary of the scatterer and $m=51$ source points on an exterior circle with radius $1.5$. The independent initial guesses were chosen as $3\pm0.8\mathrm{i}$ and $4+0.8\mathrm{i}$ for which the corresponding output is shown in \Cref{fig7} within the ranges $n\in[4,32]$ and $n\in[4,20]$, respectively. We can still recognize a one-to-one correspondence between DELs and complex-conjugated pairs of complex-valued ITE trajectories as before.

\begin{figure}[htbp]
	\centering
	\includegraphics[width=10cm]{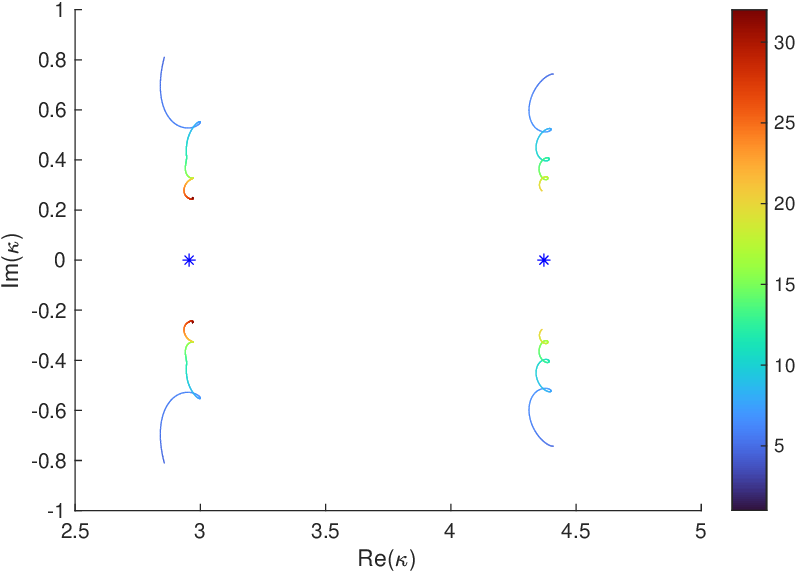}
	\caption{\label{fig7}The first two complex-conjugated pairs of  complex-valued ITE trajectories for the deformed ellipse from \Cref{deform} 
	using $n\in [4,32]$ and $n\in [4,20]$, respectively.}
\end{figure}

\subsection{Clover}
Next, we examine a clover parameterized for $t\in [0,2\pi)$ by
\begin{align}\label{clover}
	t \mapsto 0.25(\sin(4t)+5)
	\begin{pmatrix}
		\cos(t)\\
		\sin(t)
	\end{pmatrix}\ .
\end{align}
whose concrete shape is depicted in \Cref{tikz2}. 
\begin{figure}[htbp]
	\begin{center}
		\begin{tikzpicture}
			\begin{axis}[scale=0.5,
				trig format plots=rad,
				axis equal,
				hide axis,
				legend columns=4, 
				legend style={
					legend image post style={thick},
					nodes={scale=0.75},
					at={(0.5,0.2)},anchor=south,
				},
				]
				\addplot [domain=0:2*pi, samples=200, gray,opacity=1.0,fill=gray,fill opacity=0.4] ({0.25*(sin(4*x)+5)*cos(x)}, {0.25*(sin(4*x)+5)*sin(x)});
				\addplot [domain=0:2*pi, samples=200, black] ({0.25*(sin(4*x)+5)*cos(x)}, {0.25*(sin(4*x)+5)*sin(x)});
			\end{axis}
		\end{tikzpicture}
	\end{center}
	\caption{The clover parameterized by \Cref{clover}
	}\label{tikz2}
\end{figure}
We use again the modified MFS both for computing the first two DEL as well as for generating complex-valued ITE trajectories with $m_I=20$ interior points distributed on an inner circle with radius $0.2$, $m=51$ boundary collocation points and the same number of exterior source points along \Cref{clover} but scaled with a factor of $1.2$. The independent initial ITE guesses at $n=4$ were chosen as $2\pm 0.5\mathrm{i}$ and $3\pm 0.5\mathrm{i}$, respectively. The output is shown in \Cref{fig8}.
The DELs are approximately given by $2.0707$ and $3.2219$, where the second one is of particular importance as it has a numerical eigenvalue multiplicity of 2. Again, we can find two complex-conjugated pairs of complex-valued ITE trajecories approaching one of the determined DELs each, cf. \Cref{big2}, but additionally note that the eigenvalue multiplicity of each ITE from the first trajectory is one and from the second is two, resulting in overlapping for the latter. This accompanying multiplicity observation strengthens the former one-to-one correspondence between DELs and complex-conjugated pairs of complex-valued ITE trajectories even further.

\begin{figure}[htbp]
	\centering
	\includegraphics[width=10cm]{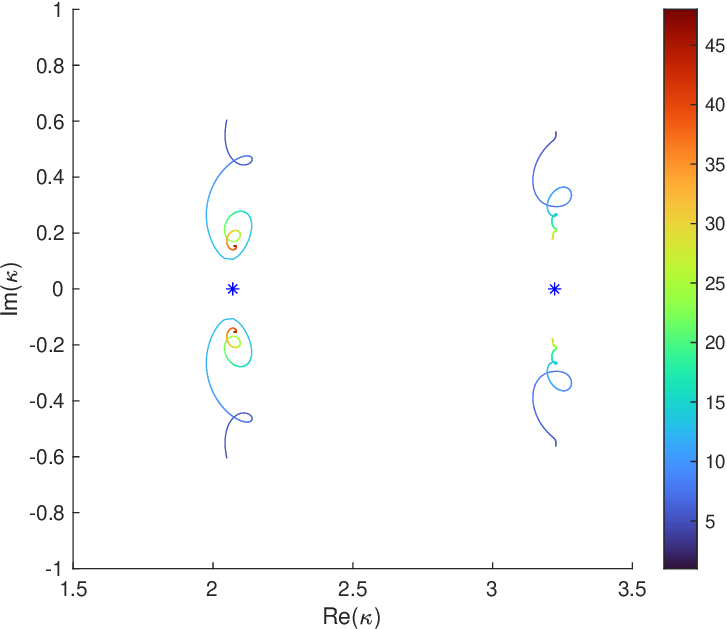}
	\caption{The first two complex-conjugated pairs of  complex-valued ITE trajectories for the clover from \Cref{clover} 
	using $n\in [4,48]$ and $n\in [4,32]$, respectively}\label{fig8}
\end{figure}

\subsection{Ellipsoid}
\begin{figure}[!ht]
	\centering
	\includegraphics[width=10cm]{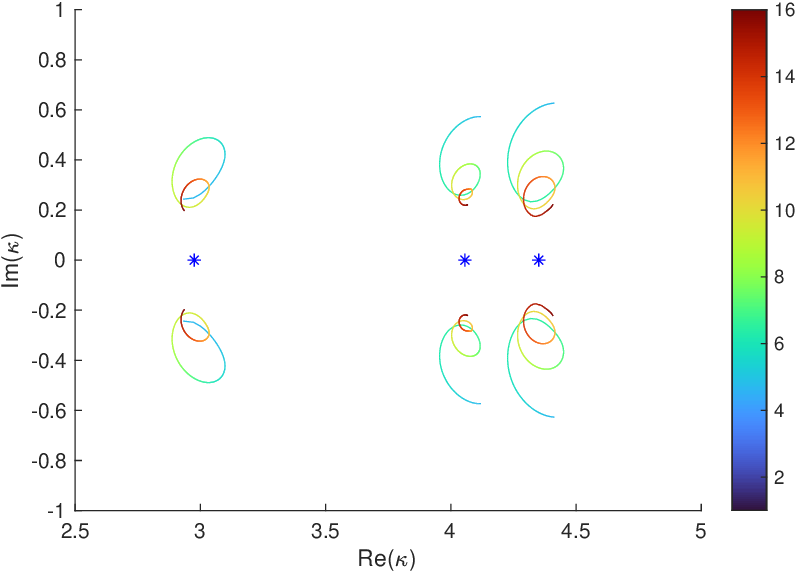}
	\caption{The first three complex-conjugated pairs of  complex-valued ITE trajectories for the ellipsoid with semi-axis $1$, $1$, and $1.2$ using $n\in [4,16]$}\label{fig11}
\end{figure}

Finally, we consider an ellipsoid in 3D with semi-axes $1$, $1$, and $1.2$. We use $m_I=41$ interior points on an inner sphere with radius $0.5$ and with equally-distributed longitude and latitude angles. Likewise, we selected $m=221$ boundary collocation points. The exterior source points are obtained by scaling the boundary collocation points with a factor of $3$. In \Cref{fig11}, the three complex-conjugated pairs of complex-valued ITE trajectories are shown for $n\in [4,16]$ and the initial starting guesses $3+0.5\mathrm{i}$, $4+0.5\mathrm{i}$, and $5+0.5\mathrm{i}$ for $n=4$, respectively.

Again, we observe a spiral behavior of the curves towards the first three IDEs, cf. \Cref{big2}, which are given as $2.9755$, $4.0563$, and $4.3511$ (cf. \cite[Table 11]{kleefeld2013numerical} for the values and their multiplicity). The first two DELs are simple whereas the third DEL has a numerical multiplicity of two. As for the clover, we observe a corresponding multiplicity behavior for the ITE trajectories which underlines the one-to-one correspondence between DELs and complex-conjugated pairs of complex-valued ITE trajectories including multiplicity also in this three-dimensional case. We feel at this point like we have collected enough samples of simply-connected scatterers all of which admit the same characteristics for DELs and complex-valued ITE trajectories to finally formulate a more general conjecture below.


\section*{Conclusion}
We have introduced the concept of ITE trajectories for the ITP subject to homogeneous media. For the unit disc as scatterer, we have proven that there is a one-to-one correspondence including multiplicity between complex-valued ITE trajectories and DELs. For more general scatterers, we could show that the only accumulation points of complex-valued ITE trajectories are restricted to DELs, too. Our numerical results for simply-connected scatterers even indicate that the one-to-one correspondence including geometric multiplicity might hold as well. Hence, we conjecture the following:

\begin{conj} 
	There is a one-to-one correspondence between complex-valued ITE trajectories and DELs of $D$ whenever $D$ is a bounded and simply-connected domain. More precisely, any complex-valued ITE trajectory $\kappa_n$ converges to some DEL as $n\to\infty$. Conversely, for any DEL there exist exactly as many complex-conjugated pairs of complex-valued ITE trajectories as the DEL's geometric multiplicity and all converge to that DEL as $n\to\infty$.
\end{conj}
Its further investigation will be subject of future research. For example, our restriction to simply-connected scatterers is due to an observation for annulus-shaped scatterers whose ITE trajectories seem to be only piece-wise complex-valued according to numerical studies. From a technical point of view, we have seen that the modified method of fundamental solution produces inaccurate results for large wave numbers and in particular for large $n$. Hence, we will investigate on how to circumvent this problem probably through the use of the computationally more demanding boundary element collocation method. 
Finally, the electromagnetic or elastic interior transmission problem will be investigated as well as a generalization of the current results to inhomogeneous media.

\section*{Acknowledgments}
The authors thank Andreas Kirsch for his helpful suggestions and fruitful discussions in the final stages of this work.
Lukas Pieronek was funded by the Deutsche Forschungsgemeinschaft (DFG, German Research Foundation) – Project-ID 258734477 – SFB 1173.

\section*{References}	
\bibliographystyle{abbrv}
\bibliography{references}

\end{document}